\documentclass[12pt,a4paper]{amsart}
\usepackage[colorlinks, linkcolor=blue, anchorcolor=blue, citecolor=green]{hyperref}
 \usepackage[normalem]{ulem}
\usepackage{graphics,epic}
\usepackage{amsmath,amssymb, amsthm,mathrsfs}
\usepackage[all,2cell]{xy}
\usepackage{multicol}
\usepackage{multirow}
\usepackage{caption}
\usepackage{subcaption}
\usepackage{shorttoc}
\usepackage{url}
\usepackage{tikz}
\usetikzlibrary{matrix,positioning,decorations.markings,arrows,decorations.pathmorphing,	
	backgrounds,fit,positioning,shapes.symbols,chains,shadings,fadings,calc}
\tikzset{->-/.style={decoration={  markings,  mark=at position #1 with
			{\arrow{>}}},postaction={decorate}}}
\tikzset{-<-/.style={decoration={  markings,  mark=at position #1 with
			{\arrow{<}}},postaction={decorate}}}

\newtheorem{theorem}{Theorem}[section]
\newtheorem*{theorem*}{Theorem}

\newtheorem{lemma}[theorem]{Lemma}

\newtheorem*{conjecture*}{Conjecture}

\newtheorem{example}[theorem]{Example}
\newtheorem{remark}[theorem]{Remark}
\newtheorem{definition}[theorem]{Definition}

\newtheorem{thm}[theorem]{Theorem}
\newtheorem{lem}[theorem]{Lemma}
\newtheorem{prop}[theorem]{Proposition}
\newtheorem{cor}[theorem]{Corollary}
\newtheorem{conj}[theorem]{Conjecture}

\newcommand{\ie}{{\em i.e.}\ }

\newcommand{\opname}[1]{\operatorname{\mathsf{#1}}}

\renewcommand{\mod}{\opname{mod}\nolimits}

\newcommand{\D}{\mathbb{D}}

\renewcommand{\P}{\mathbb{P}}

\newcommand{\fM}{{\mathbf M}}
\newcommand{\fS}{{\mathbf S}}
\newcommand{\fP}{{\mathbf P}}
\newcommand{\fT}{{\mathbf T}}

\newcommand{\Arcseg}{\mathbf{Arcs}}
\newcommand{\Seg}{\mathbf{Seg}}
\newcommand{\ba}{\mathbf{a}}
\newcommand{\Int}{\mathbf{Int}}
\newcommand{\Intv}{\mathbf{\underline{Int}}}
\newcommand{\dimv}{\mathbf{\underline{dim}}}

\newcommand{\TT}{\mathbf{T}}
\renewcommand{\SS}{\mathbf{S}}
\newcommand{\MM}{\mathbf{M}}

\newcommand{\Hom}{\opname{Hom}}

\newcommand{\Ext}{\opname{Ext}}

\newcommand{\End}{\opname{End}}

\newcommand{\dInt}{\opname{Int}}

\newcommand{\dInttv}{\underline{\opname{Int}}}

\newcommand{\cc}{{\mathcal C}}

\newcommand{\cm}{{\mathcal M}}
\newcommand{\cM}{{\mathcal M}}
\newcommand{\cn}{{\mathcal N}}

\newcommand{\E}{{\mathbb E}}

\newcommand{\mc}{\mathcal{C}}

\newcommand{\ma}{\mathcal{A}}

\newcommand{\Arc}{\mathbf{Arcs}}
\newcommand{\Arcs}{\mathbf{Arcs}}

\newcommand{\nn}{node[black]{$\bullet$}}

\makeatletter
\newcommand*\bigcdot{\mathpalette\bigcdot@{.5}}
\newcommand*\bigcdot@[2]{\mathbin{\vcenter{\hbox{\scalebox{#2}{$\m@th#1\bullet$}}}}}
\makeatother

\begin{document}

\title{On denominator conjecture for cluster algebras of finite type}

\author{Changjian Fu}
\address{Changjian Fu\\Department of Mathematics\\SiChuan University\\610064 Chengdu\\P.R.China}
\email{changjianfu@scu.edu.cn}
\author{Shengfei Geng}
\address{Shengfei Geng\\Department of Mathematics\\SiChuan University\\610064 Chengdu\\P.R.China}
\email{genshengfei@scu.edu.cn}
\keywords{denominator conjecture, finite type, marked surface}

\begin{abstract}
We continue our investigation on denominator conjecture of Fomin and Zelevinsky for cluster algebras via geometric models initialed in \cite{FG22}. 
In this paper, we confirm the denominator conjecture for cluster algebras of finite type. 
Our main contribution is a proof of this conjecture for cluster algebras of type $\mathbb{D}$, along with an algorithm for the exceptional types. For the type $\mathbb{D}$ cases, our approach involves a geometric model based on discs with a puncture. By removing the puncture or changing it to an unmarked boundary component, we also provide an alternative proof of the denominator conjecture for cluster algebras of types $\mathbb{A}$ and $\mathbb{C}$, respectively.

\end{abstract}

\maketitle
\tableofcontents

\section{Introduction}
A cluster algebra \cite{FZCA1} is a commutative algebra with a distinguished set of generators called {\em cluster variables}. These generators are gathered into overlapping sets of fixed finite cardinality, called {\em clusters}, which are defined recursively from an initial one via {\em mutations}. 
A monomial in cluster variables belonging to the same cluster is called a {\em cluster monomial}. Let $\mathcal{A}$  be a cluster algebra with trivial coefficients, and let  $\mathbf{x}=(x_1,\dots, x_n)$ be an initial cluster. According to the Laurent phenomenon \cite{FZCA1}, every cluster monomial $z$ can be uniquely expressed as
\[
z=\frac{f(x_1,\dots, x_n)}{x_1^{d_1}\cdots x_n^{d_n}},
\]
where $d_1,\dots, d_n\in \mathbb{Z}$ and $f(x_1,\dots, x_n)\in \mathbb{Z}[x_1,\dots, x_n]$ is not divisible by any $x_i$. The integer vector $d_\mathbf{x}(z)=(d_1,\dots, d_n)$ is called the {\em denominator vector} of $z$ with respect to $\mathbf{x}$. We remark that the denominator vector $d_{\mathbf{x}}(z)$ depends on the choice of the initial cluster $\mathbf{x}$.

Inspired by Lusztig’s parameterization of canonical bases in the theory of quantum groups, Fomin and Zelevinsky \cite{FZ03,FZCA4} formulated the following denominator conjecture, which remains a challenging problem in the field of cluster algebras.
\begin{conj}\cite[Conjecture 4.17]{FZ03}\cite[Conjecture 7.6]{FZCA4}\label{conj:denominator-conj}
   Different cluster monomials have different denominator vectors with respect to any given cluster.
\end{conj}
Conjecture \ref{conj:denominator-conj} was verified for cluster algebras of rank $2$ (\ie $n=2$) via combinatorics of root systems of corresponding Kac-Moody algebras \cite{SZ04} and acyclic cluster algebras with respect to an acyclic initial cluster by representation theory of (valued) quivers \cite{CK06,CK08,RS20}.
Based on geometric models provided by marked surfaces,  
  Conjecture \ref{conj:denominator-conj} was established for cluster algebras of type $\mathbb{A}$, $\mathbb{B}$ and $\mathbb{C}$ in \cite{FG22}, and a large class of cluster algebras of surface type with respect to particular choices of initial clusters in \cite{FG24}. A weaker version of this conjecture has also been explored in several works, see \cite{GP12,NS14,FuG19,GyYu20,FGL21} for instance.  In contrast, Jiarui Fei \cite{Fei} found a counterexample to Conjecture \ref{conj:denominator-conj}. The associated Jacobian algebra of Fei's counterexample is of wild type. Nevertheless, one may still expect that this conjecture is true for many important classes of cluster algebras.

The aim of this note is to establish Conjecture \ref{conj:denominator-conj} for cluster algebras of finite type.
\begin{thm}\label{thm:denominator-conj-finite-type}
    Let $\mathcal{A}$ be a cluster algebra of finite type. Then different cluster monomials have different denominator vectors with respect to any given cluster.
\end{thm}
Recall that a cluster algebra is of {\em finite type} if it has finitely many cluster variables. Fomin and Zelevinsky \cite{FZCA2} proved that cluster algebras of finite type are parameterized by finite root systems.  According to \cite{FG22} and \cite{SZ04}, it remains to verify Conjecture \ref{conj:denominator-conj} for cluster algebras of type $\mathbb{D}$ and exceptional type $\mathbb{E}_6$, $\mathbb{E}_7$, $\mathbb{E}_8$ and $\mathbb{F}_4$.  For cluster algebras of type $\mathbb{D}$, we use geometric models based on discs with a puncture. By the Fomin--Shapiro--Thurston correspondence \cite{FST08}, there is a bijection between the set of cluster monomials and the set of finite multisets of pairwise compatible admissible tagged arcs.
Under this correspondence, the denominator vector of a cluster monomial corresponds to the intersection vector of the corresponding finite multiset. We prove in Theorem \ref{t: main theorem} that such a multiset is uniquely determined by its intersection vector. By removing the puncture or changing the puncture to an unmarked boundary component, this also provides an alternative, simple proof of Conjecture \ref{conj:denominator-conj}  for type $\mathbb{A}$ and $\mathbb{C}$. For the exceptional types, we verify Conjecture \ref{conj:denominator-conj} by using an algorithm.

Note that acyclic cluster algebras are categorified by cluster categories. The following is a representation-theoretic version of Theorem \ref{thm:denominator-conj-finite-type}, which improves upon the corresponding results in \cite{GP12,R11}.
\begin{cor}\label{cor:cluster-tilted-alg}
    Let $\Gamma$ be a cluster-tilted algebra of finite representation type. Different $\tau$-rigid $\Gamma$-modules have different dimension vectors.
\end{cor}

When preparing this article, Yurikusa kindly shared with us his work \cite{Yu}, in which he also verified Conjecture \ref{conj:denominator-conj} for a large class of cluster algebras of surface type with respect to particular choices of initial clusters, using a different method. In particular, he also obtained that Conjecture \ref{conj:denominator-conj} holds for cluster algebras of type $\mathbb{D}$.

The paper is organized as follows. In Section \ref{s:disc}, we recollect definitions related to marked surfaces. Since we are only interested in cluster algebras of finite type, we restrict ourselves to discs with at most one puncture.  Section \ref{s:local-arc-segmet} is devoted to showing that the intersection vector of a finite permissible multiset $\cm$ uniquely determines the multiset of irreducible arc segments associated with $\cm$ (cf. Theorem \ref{thm:intersection-arc-segments}). By applying Theorem \ref{thm:intersection-arc-segments}, we prove in Section \ref{s:inter-vector-determine-arc} that a finite multiset is uniquely determined by its intersection vector (Theorem \ref{t: main theorem}). A similar result also holds for a disc with an unmarked boundary component in its interior (Proposition \ref{thm:type-c-intersection-vector}). After recalling the Fomin--Shapiro--Thurston correspondence in Section \ref{ss:FST-corr}, we prove Theorem \ref{thm:denominator-conj-finite-type} in Section \ref{ss:proof-main-result} and Corollary \ref{cor:cluster-tilted-alg} in Section \ref{ss: cluster-tilted-algebra}.

\subsection*{Acknowledgements}
Both authors are grateful to Professor Siyang Liu for informing us Fei's counterexample, to  Professor Jiarui Fei for explaining us the counterexample and to Professor Toshiya Yurikusa for sharing us his preprint \cite{Yu}. We also thank Professor Jianrong Li and Professor Xuexue Zhou for helpful discussion on the program. This work is partially support by the National Natural Science Foundation of China (Grant No. 11971326, 12471037).

\section{Disc with at most one puncture}\label{s:disc}
\subsection{Disc with at most one puncture}\label{ss:disc-one-puncture}
We follow \cite{FST08}.
Let $\fS$ be a disc, $\fM\subset \fS$ a non-empty finite set of marked points on the boundary $\partial \fS$ and $\fP\subset \fS\backslash \partial \fS$  a finite set
 of punctures in the interior of $\fS$. 
We assume that $|\fM|\geq 4$ and $|\fP|\leq 1$. When $|\fP|=0$, then $\fS:=(\fS,\fM)$ is called an {\it $n$-gon}, where $|\fM|=n$. When $|\fP|=1$, then $\fS:=(\fS,\fM,\fP)$ is called a {\it once-punctured disc}. In this case, we always denote by $P$ the unique puncture.


A {\em curve} on $\fS$ is a continuous map $\gamma: [0,1]\to \fS$ such that
\begin{itemize}
    \item[$\circ$] $\gamma(0),\gamma(1)\in \fM\cup \fP$ and $\gamma(t)\in \fS\backslash (\fM\cup \fP)$ for $0<t<1$;
    \item[$\circ$] $\gamma$ is neither null-homotopic nor homotopic to a boundary segment.  
\end{itemize}
We always consider curves on $\fS$ up to homotopy relative to their endpoints.
The inverse of a curve $\gamma$ on $\fS$ is defined as $\gamma^{-1}(t)=\gamma(1-t)$ for $t\in [0,1]$. A curve $\gamma$ is called a {\it loop} if $\gamma(0)=\gamma(1)$.
Throughout this paper, when we consider intersections of curves, we assume that they intersect transversally in a minimum number of points in $\fS\backslash\fM$. 
For any curves $\gamma_1, \gamma_2$ on $\fS$,
the {\it intersection number} between $\gamma_1$ and $\gamma_2$ is defined to be
\[\Int^\circ(\gamma_1,\gamma_2):=|\{(t_1,t_2)~|~0<t_1,t_2<1, \gamma_1(t_1)=\gamma_2(t_2)\}|.
\]
 We say that  $\gamma_1$ and $\gamma_2$ are {\it compatible} if $\Int^\circ(\gamma_1,\gamma_2)=0$.  An {\it arc} $\gamma$ on $\fS$
 is a curve without self-intersections, that is, $\Int^\circ(\gamma,\gamma)=0$. We denote by $\mathbb{A}(\fS)$ the set of equivalence classes of  arcs on $\fS$  under the equivalence relation given by taking inverse.

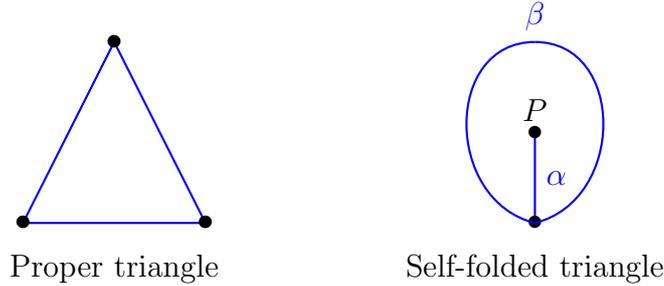
\begin{figure}[ht]
\begin{minipage}[t]{0.3\linewidth} 

\begin{tikzpicture}[xscale=0.6,yscale=0.6]

			\draw[blue,thick] (-2,-2)\nn to (2,-2)\nn;
			
		\draw[blue,thick] (-2,-2)\nn to (0,2)\nn;
			\draw[blue,thick] (2,-2)\nn to (0,2)\nn;

\node at (0,-3){Proper triangle};
	\end{tikzpicture}
\end{minipage}%
\begin{minipage}[t]{0.3\linewidth} 
\centering
\begin{tikzpicture}[xscale=0.6,yscale=0.6]
			
			\draw[blue,thick] (0,-2)\nn to (0,0)\nn;
			
\draw[thick,color=blue] (0,-2) .. controls (-2,-1.5) and (-2,2) .. (0,2);	
\draw[thick,color=blue] (0,-2) .. controls (2,-1.5) and (2,2) .. (0,2);		
\draw[thick,blue](0,-1)node[right]{$\alpha$}(0,2)node[above]{$\beta$};
\node at (0,.5){$P$};

\node at (0,-3){Self-folded triangle};

	\end{tikzpicture}
\end{minipage}

\caption{Basic tiles of an ideal triangulation}\label{f:triangle}
\end{figure}

 An {\it ideal triangulation} $\fT$ of $(\fS,\fM)$ is a maximal set of pairwise compatible arcs.  An ideal triangulation $\fT$ of $\fS$ decomposes $\fS$
 into triangles.
 A triangle in $\fT$ has either three distinct sides  or  it is a self-folded triangle, cf.  
 Figure~\ref{f:triangle}, where $P$ is the puncture, and we call $\alpha$ the {\em folded side} and $\beta$ the {\em remaining side}.
In particular, if $\fS$ is an $n$-gon, then an ideal triangulation decomposes $\fS$ into  proper triangles.

\subsection{Tagged arc and tagged triangulation}
Each arc $\gamma$ in $(\fS, \fM)$ has two ends obtained by arbitrarily cutting $\gamma$ into three pieces, and then throwing out the middle one. We think of the two ends as locations near the endpoints to be used for labeling (tagging) an arc. 

A {\em tagged arc} is an arc in which each end has been tagged in one of two ways, {\em plain} or {\em notched}, so that the following conditions are satisfied:
 \begin{itemize} 
\item an endpoint lying on the boundary is tagged plain; 
\item both ends of a loop are tagged in the same way.
\end{itemize}
We denote by $\mathbb{A}_{\bowtie}(\fS)$ the set of tagged arcs on $\fS$. 
In the figures, we represent tags as follows:
\begin{figure}[ht]
\begin{minipage}[t]{0.3\linewidth} 
\begin{tikzpicture}[xscale=0.6,yscale=0.6]
\draw[blue,thick] (-4,0) to (-2,0)\nn;
\draw[blue,thick](2,0) to (4,0)\nn;
\node at (-5,0){plain};
\node at (0.5,0){notched};
\node at (3.5,0)[rotate=90]{$\bowtie$};
\end{tikzpicture}
\end{minipage}
\end{figure}

Let $\gamma\in \mathbb{A}_{\bowtie}(\fS)$,
 \begin{itemize}
\item $\gamma$ is called a {\it plain arc} if both ends are tagged plain;
\item $\gamma$ is called a {\it 1-notched arc} if  one end of $\gamma$ is tagged plain and the other end is tagged notched.
\end{itemize}

For  $\alpha\in \mathbb{A}_{\bowtie}(\fS)$, we denote by $\bar{\alpha}$ its underlying untagged arc.
For  $\alpha,\beta\in \mathbb{A}_{\bowtie}(\fS)$ such
that $\bar{\alpha}=\bar{\beta}$, if exactly one of them is  $1$-notched, then the pair $(\alpha,\beta)$ is called a {\it  pair of conjugate arcs}.  Since  we have at most one puncture $P$, each pair of conjugate arcs has different tags at $P$. If the other endpoint is $Q$, we may denote such a pair by $(\gamma_Q^-,\gamma_Q^{\bowtie})$, 
where $\gamma_Q^-$ is tagged plain at the end incident to  $P$, while $\gamma_Q^{\bowtie}$ is tagged notched at the end incident to $P$.

Two tagged arcs $\alpha,\beta\in \mathbb{A}_{\bowtie}(\fS)$ are called {\it compatible} if and only if the following conditions are satisfied:
\begin{itemize}
\item  $\bar{\alpha}$ and $\bar{\beta}$ are compatible, \ie $\Int^\circ(\bar{\alpha}|\bar{\beta})=0$;
\item if  $\bar{\alpha}\neq \bar{\beta}$ and $\alpha, \beta$ share an endpoint $Q$, then the ends of $\alpha$ and $\beta$ connected to $Q$ must be tagged in the same way;
\item if $\bar{\alpha}=\bar{\beta}$, then at least one end of $\alpha$ must be tagged in the same way as the corresponding end of $\beta$.
\end{itemize}

A tagged arc that does not cut out a once-punctured monogon is called {\it admissible}.
A maximal collection $\fT$ of pairwise compatible admissible tagged arcs  on $\fS$ is
called a {\em tagged triangulation.}  Since $\fS$ is a disk with at most one puncture and each tagged arc in $\fT$ is admissible, $\fT$ contains no loops.
A tagged triangulation $\fT$ of $\fS$ decomposes $\fS$ into  two types of basic tiles as shown in Figure \ref{f:basic tile}.
\begin{figure}[ht]
\begin{minipage}[t]{0.4\linewidth} 

\begin{tikzpicture}[xscale=0.6,yscale=0.6]

		\draw[blue,thick] (-2,-2)\nn to (2,-2)\nn;
			
		\draw[blue,thick] (-2,-2)\nn to (0,2)\nn;
			\draw[blue,thick] (2,-2)\nn to (0,2)\nn;
			
		\draw[blue](-2,-2)node[left]{$P_2$}(2,-2)node[right]{$P_3$}(0,2)node[above]{$P_1$};

\node at (0,-4){Type I:  triangle with };
\node at (0,-5){three different vertices};
	\end{tikzpicture}
\end{minipage}%
\begin{minipage}[t]{0.4\linewidth} 
\centering
\begin{tikzpicture}[xscale=0.6,yscale=0.6]

\draw[blue,thick] (0,2)\nn to (0,2)\nn;
			\draw[blue,thick] (0,-2)\nn to (0,-2)\nn;
				\draw[blue,thick] (0,0)\nn to (0,0)\nn;
\draw[thick,color=blue] (0,-2) .. controls (-2,-1) and (-2,1) .. (0,2);	
\draw[thick,color=blue] (0,-2) .. controls (2,-1) and (2,1) .. (0,2);	

\draw[thick,color=blue] (0,-2) .. controls (-0.3,-1.4) and (-0.3,-0.6) .. (0,0);	
\draw[thick,color=blue] (0,-2) .. controls (0.3,-1.4) and (0.3,-0.6) .. (0,0);

		
\draw[blue](0,-2)node[below]{$P_1$}(0,2)node[above]{$P_2$}(0,0)node[right]{$P$};	
\draw[blue](0.17,-0.4)node[rotate=20]{$\bowtie$};			

\node at (0,-4){Type II: 1-puncture piece, where};
\node at (0,-5){$P_1$ and $P_2$ are different marked points};

	\end{tikzpicture}
\end{minipage}

\caption{Basic tiles of a tagged triangulation}\label{f:basic tile}
\end{figure}

\begin{remark}
    If $\fS$ is an $n$-gon, then a tagged arc on $\fS$ is the same as an arc, and a tagged triangulation is the same as an ideal triangulation. In particular, only type I tiles will appear.
\end{remark}

\subsection{Intersection number and Intersection vector}
The following definition plays a central role in our investigation.
\begin{definition}\cite[Definition 8.4]{FST08}\label{def:intersection-number-tagged-arc}
Let $\alpha$ and $\beta$ be two tagged arcs on $\fS$. 
The {\em intersection number} $\Int(\alpha|\beta)$ of $\alpha$ and $\beta$ is defined by
\[\Int(\alpha|\beta)= \Int^A(\alpha|\beta)+\Int^B(\alpha|\beta)+\Int^C(\alpha|\beta)+\Int^D(\alpha|\beta),\]
where
\begin{itemize}
\item $\Int^A(\alpha|\beta)=\Int^\circ(\bar{\alpha}|\bar{\beta})$;

\item $\Int^B(\alpha|\beta)=0$ unless $\bar{\alpha}$ is a loop based at a point $a\in \fM\cup \fP$, in which case $\Int^B(\alpha|\beta)$ is computed
as follows: assume that $\bar{\beta}$ intersects $\bar{\alpha}$ {\rm(}in the interior of $\fS\backslash (\fM\cup \fP)${\rm)} at the points $b_1,\ldots,b_m$ {\rm(}numbered along $\bar{\beta}$ in this order{\rm )}; then $\Int^B(\alpha|\beta)$ is the negative of the number of segments $\gamma_i=[b_i, b_{i+1}]\subset \beta$ such that $\gamma_i$ together with the segments $[a, b_i]\subset \alpha$ and $[a, b_{i+1}]\subset \alpha$ form three sides of a contractible triangle disjoint from the punctures; 
\item $\Int^C(\alpha|\beta)=0$ unless $\bar{\alpha}=\bar{\beta}$, in which case $\Int^C(\alpha|\beta)=-1$. 

\item $\Int^D(\alpha|\beta)$ is the number of ends of $\beta$ that are incident to an endpoint of $\alpha$ and, at that endpoint, carry a tag different from the one $\alpha$ does. 
\end{itemize}
\end{definition}
\begin{remark}
It is easy to see that $\alpha,\beta\in \mathbb{A}_{\bowtie}(\fS)$ are compatible if and only if $\Int(\alpha|\beta)\leq 0$.  Moreover, when $\alpha\neq \beta$, $\alpha$ and $\beta$ are compatible if and only if $\Int(\alpha|\beta)=0$.
\end{remark}

\begin{example}
Assume that $\fS$ is an $n$-gon and $\alpha,\beta\in \mathbb{A}(\fS)
$. We have $\Int^B(\alpha|\beta)=0=\Int^D(\alpha|\beta)$.
It follows that $\Int(\alpha|\beta)=\Int^A(\alpha|\beta)+\Int^C(\alpha|\beta)$.
\end{example}

Let $\fT$ be a tagged triangulation on $\fS$ and $\alpha\in \mathbb{A}_{\bowtie}(\fS)$. The {\em intersection vector} $\underline{\Int}_\fT(\alpha)$ of $\alpha$ with respect to $\fT$ is defined as
\begin{equation}\label{eq:inter-vector-arc}
    \underline{\Int}_\fT(\alpha)=(\Int(\mathbf{a}|\alpha))_{\mathbf{a}\in \fT}.
\end{equation}
For a finite multiset $\cM$ of pairwise compatible tagged arcs, we define
\begin{eqnarray}
    \Int(\mathbf{a}|\cM)&:=&\sum_{\alpha\in \cM}\Int(\mathbf{a}|\alpha), \mathbf{a}\in \fT\label{eq:inter-number-multiset} \\
    \underline{\Int}_\fT(\cM)&:=&\sum_{\alpha\in \cM}\underline{\Int}_\fT(\alpha). \label{eq:inter-vector-multiset}
\end{eqnarray}
By definition, we have $\underline{\Int}_\fT(\cM)=(\Int(\mathbf{a}|\cM))_{\mathbf{a}\in \fT}$.

\subsection{Permissible  multiset}\label{ss: permissible-multiset}
Let $\fT$ be a tagged triangulation of $(\SS,\MM)$.
A finite multiset $\mathcal{C}$ of tagged arcs is  {\em permissible } with respect to $\fT$ if $\mathcal{C}$ satisfies the following conditions:
\begin{itemize}
\item[($\mathbf{P1}$)] Any two tagged arcs in $\mathcal{C}$ are compatible with each other;
\item[($\mathbf{P2}$)] $\mathcal{C}$  contains no conjugate pairs of tagged arcs;
\item[($\mathbf{P3}$)]  $\mathcal{C}$ contains no once-punctured loop that closely wraps  around an arc in $\fT$;
\item[($\mathbf{P4}$)] $\mathcal{C}\cap \fT=\emptyset$.
\end{itemize}
A  tagged arc $\gamma$ is permissible if the set $\{\gamma\}$ is permissible with respect to $\fT$.

For a finite multiset $\cm$ consisting of pairwise compatible admissible tagged arcs, denote by $\cm^{*}$ the new multiset obtained by replacing each conjugate pair $(\gamma_Q^-,\gamma_Q^{\bowtie})$ in $\cm$ with the corresponding loop $l_{Q}$ based at $Q$   such that $l_Q$ closely  wraps  around $\gamma_Q^-$.  If $\cm\cap \fT=\emptyset$, then
one can deduce that  $\cm^{*}$ is a permissible  multiset  of tagged arcs with respect to $\fT$. Assume that  $\cm$ contains at least a pair of conjugate tagged arcs  $(\gamma_Q^-,\gamma_Q^{\bowtie})$,   since all tagged arcs in $\cm$ are  pairwise compatible  with each other,  one can deduce that any other tagged arc in $\cm$  incident to $P$  is either a copy of $\gamma_Q^-$ or a copy of $\gamma_Q^{\bowtie}$.


\begin{lemma}\label{l: m* has same intersection vector with m}
Let  $\cm$ be a finite multiset  consisting of pairwise compatible admissible tagged arcs such that $\cm\cap \fT=\emptyset$.
Then  $\Intv_{\fT}(\cm)=\Intv_{\fT}(\cm^{*})$.
\end{lemma}
\begin{proof}
Without loss of generality, we may assume that $\cm$ contains at least a pair of conjugate tagged arcs.  Otherwise, $\cm=\cm^\ast$ and the result follows. 

Let $(\gamma_Q^-,\gamma_Q^{\bowtie})$ be a pair of conjugate tagged arcs with $Q$ as their boundary endpoint, 
$l_{Q}$ the loop  which closely  wraps around $\gamma_Q^-$. By the construction of $\cm^*$,
it suffices to prove that for each $\ba\in\fT$, $$\Int(\ba|\gamma_Q^-)+\Int(\ba|\gamma_Q^{\bowtie})=\Int(\ba|l_Q).$$
Since $\cm\cap \fT=\emptyset$,
we know that $\cm^*\cap \fT=\emptyset$.  Note that $\fT$ contains no loops, we have
\begin{align}\label{gongshi:intersection of a and b}
  \Int(\ba|\beta)=\Int^A(\ba|\beta)+\Int^D(\ba|\beta)  
\end{align}
for any tagged arc $\ba\in \TT$ and $\beta$ in $\cm$ or $\cm^\ast$.
Furthermore, if $\ba$ is a tagged arc with  both endpoints  on the boundary, then $\Int(\ba|\beta)=\Int^A(\ba|\beta)$.

Assume that $\ba$ is a tagged arc with  both endpoints  on the boundary. Because $l_Q$ is a loop closely wrapping $\gamma_Q^-$, we obtain 
\begin{align*}
    \Int(\ba|l_Q)&=\Int^A(\ba|l_Q)\\
    &=2\Int^A(\ba|\gamma_Q^-)\\
    &=\Int^A(\ba|\gamma_Q^-)+\Int^{A}(\ba|\gamma_Q^{\bowtie})\\
    &=\Int(\ba|\gamma_Q^-)+\Int(\ba|\gamma_Q^{\bowtie}).
\end{align*}

Now assume that $\ba$ is a tagged arc with  the puncture $P$ as an endpoint. Because $P$ is also the endpoint of $\gamma_Q^-$ and $\gamma_Q^{\bowtie}$, we have $$\Int(\ba|\gamma_Q^{-})+\Int(\ba|\gamma_Q^{\bowtie})=\Int^D(\ba|\gamma_Q^{-})+\Int^D(\ba|\gamma_Q^{\bowtie})=1$$ 
by (\ref{gongshi:intersection of a and b}), no matter whether $\ba$ is tagged plain or notched at the end incident to $P$.
On the other hand, since $l_Q$ is a loop closely wrapping $\gamma_Q^-$ and $\ba\neq \gamma_Q^-,\ba\neq \gamma_Q^{\bowtie}$, by (\ref{gongshi:intersection of a and b}), one can also get that  $$\Int(\ba|l_Q)=\Int^A(\ba|l_Q)=1.$$ Hence  $$\Int(\ba|l_Q)=\Int(\ba|\gamma_Q^-)+\Int(\ba|\gamma_Q^{\bowtie}).$$
\end{proof}

\subsection{Irreducible arc segments}
\begin{figure}[ht]

	\begin{minipage}{0.3\linewidth} 
 \centering
		\begin{tikzpicture}[xscale=0.6,yscale=0.6]
			
			\draw[blue,thick] (-2,-1)\nn to (2,-1)\nn;
			
		\draw[blue,thick] (-2,-1)\nn to (0,3)\nn;
			\draw[blue,thick] (2,-1)\nn to (0,3)\nn;



	\draw[red,thick] (0,3) to (-0.5,-1.3);
	\node at (0,-3){Type 1.1};

		\end{tikzpicture}
		\end{minipage} 
  \begin{minipage}{0.3\linewidth} 
 \centering
		\begin{tikzpicture}[xscale=0.6,yscale=0.6]
			
			\draw[blue,thick] (-2,-1)\nn to (2,-1)\nn;
			
		\draw[blue,thick] (-2,-1)\nn to (0,3)\nn;
			\draw[blue,thick] (2,-1)\nn to (0,3)\nn;
			

	\draw[thick,color=red] (-1.1,1) .. controls (0,-0.15) and (1,1) .. (1.2,1);			
		
	\node at (0,-3){Type 1.2};

		\end{tikzpicture}
		\end{minipage} 
  \begin{minipage}{0.3\linewidth} 
 \centering
		\begin{tikzpicture}[xscale=0.6,yscale=0.6]
			
			\draw[blue,thick] (-2,-1)\nn to (2,-1)\nn;
			
		\draw[blue,thick] (-2,-1)\nn to (0,3)\nn;
			\draw[blue,thick] (2,-1)\nn to (0,3)\nn;

		\draw[blue](0,3)node[above]{$P$};	
		\draw[blue](0,0)node[left]{$\eta$}(-1,1)node[left]{$\ba$};

	\draw[red,thick] (0,3) .. controls (-0.5,0) and (-0.5,0) .. (-2,-1);
	\node at (0,-3){Type 1.3};

		\end{tikzpicture}
		\end{minipage} 
  \caption{Irreducible arc segments in a type I tile}\label{fig:type-I-irred-arc-seg}
 \end{figure}
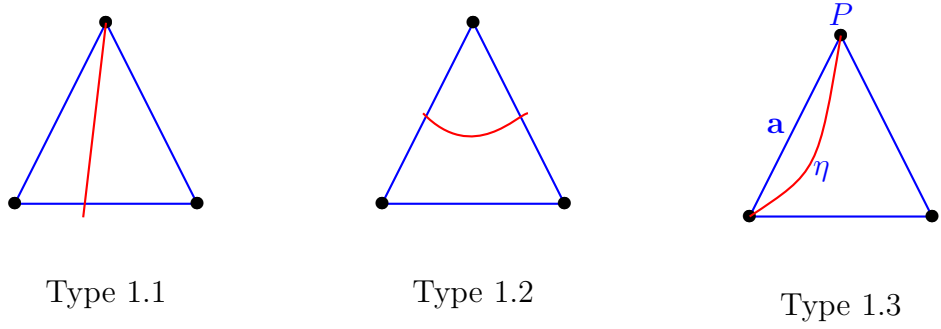

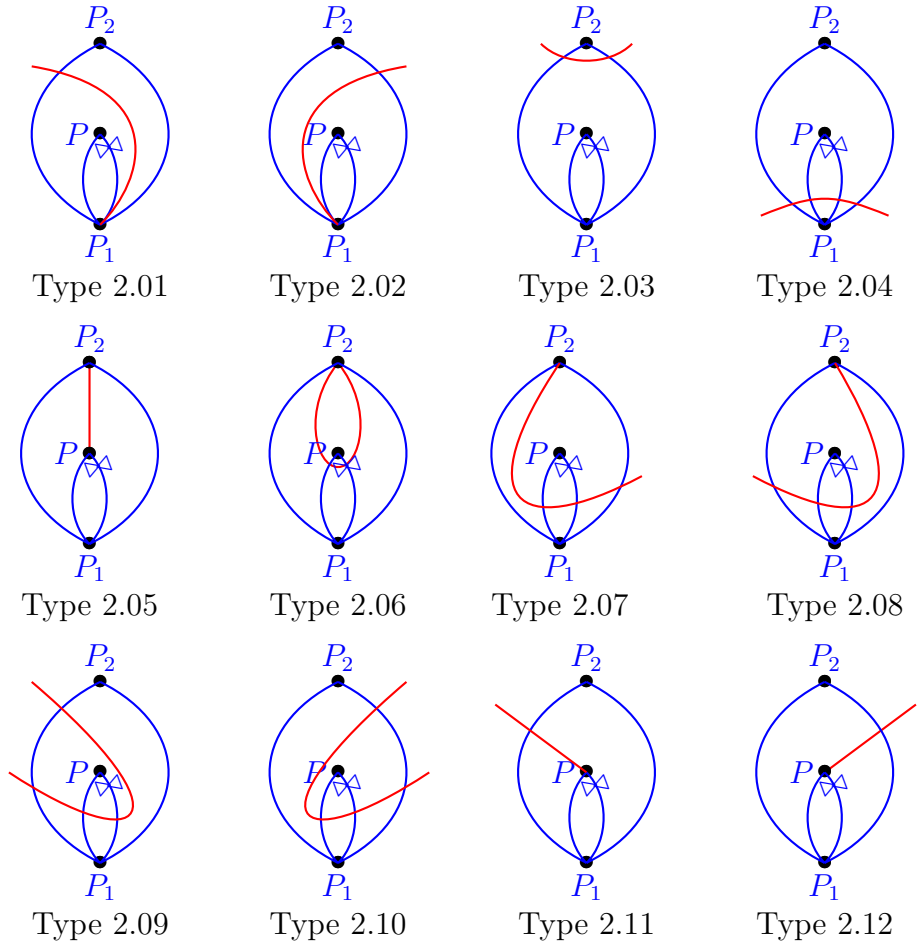
\begin{figure}[ht]    
\begin{minipage}[t]{0.2\linewidth} 
\centering
\begin{tikzpicture}[xscale=0.6,yscale=0.6]

			\draw[blue,thick] (0,-2)\nn to (0,-2)\nn;
				\draw[blue,thick] (0,0)\nn to (0,0)\nn;
				
		\draw[blue,thick] (0,2)\nn to (0,2)\nn;	
				
\draw[thick,color=blue] (0,-2) .. controls (-2,-1) and (-2,1) .. (0,2);	
\draw[thick,color=blue] (0,-2) .. controls (2,-1) and (2,1) .. (0,2);	

\draw[thick,color=blue] (0,-2) .. controls (-0.5,-1.5) and (-0.5,-0.5) .. (0,0);	
\draw[thick,color=blue] (0,-2) .. controls (0.5,-1.5) and (0.5,-0.5) .. (0,0);	
		
\draw[blue](0,-2)node[below]{$P_1$}(0,2)node[above]{$P_2$}(0,0)node[left]{$P$};	
\draw[blue](0.2,-0.3)node[rotate=20]{$\bowtie$};	

\draw[thick,color=red] (0,-2) .. controls (1,-1) and (1.5,1) .. (-1.5,1.5);
\node at (0,-3.4){Type 2.01};
	\end{tikzpicture}
	
\end{minipage}%
\begin{minipage}[t]{0.2\linewidth} 
\centering
\begin{tikzpicture}[xscale=0.6,yscale=0.6]
			\draw[blue,thick] (0,2)\nn to (0,2)\nn;
			
			\draw[blue,thick] (0,-2)\nn to (0,-2)\nn;
				\draw[blue,thick] (0,0)\nn to (0,0)\nn;

\draw[thick,color=blue] (0,-2) .. controls (-2,-1) and (-2,1) .. (0,2);	
\draw[thick,color=blue] (0,-2) .. controls (2,-1) and (2,1) .. (0,2);	

\draw[thick,color=blue] (0,-2) .. controls (-0.5,-1.5) and (-0.5,-0.5) .. (0,0);	
\draw[thick,color=blue] (0,-2) .. controls (0.5,-1.5) and (0.5,-0.5) .. (0,0);	
		
\draw[blue](0,-2)node[below]{$P_1$}(0,2)node[above]{$P_2$}(0,0)node[left]{$P$};	
\draw[blue](0.2,-0.3)node[rotate=20]{$\bowtie$};	
 \draw[thick,color=red] (0,-2) .. controls (-1,-1) and (-1.5,1) .. (1.5,1.5);

\node at (0,-3.4){Type 2.02};
		\end{tikzpicture}

\end{minipage}
\begin{minipage}[t]{0.2\linewidth} 
\centering
\begin{tikzpicture}[xscale=0.6,yscale=0.6]
			
	\draw[blue,thick] (0,2)\nn to (0,2)\nn;
			
			\draw[blue,thick] (0,-2)\nn to (0,-2)\nn;
				\draw[blue,thick] (0,0)\nn to (0,0)\nn;

\draw[thick,color=blue] (0,-2) .. controls (-2,-1) and (-2,1) .. (0,2);	
\draw[thick,color=blue] (0,-2) .. controls (2,-1) and (2,1) .. (0,2);	

\draw[thick,color=blue] (0,-2) .. controls (-0.5,-1.5) and (-0.5,-0.5) .. (0,0);	
\draw[thick,color=blue] (0,-2) .. controls (0.5,-1.5) and (0.5,-0.5) .. (0,0);	
		
\draw[blue](0,-2)node[below]{$P_1$}(0,2)node[above]{$P_2$}(0,0)node[left]{$P$};	
		
\draw[blue](0.2,-0.3)node[rotate=20]{$\bowtie$};	
	
\draw[thick,color=red] (1,2) .. controls (0.5,1.5) and (-0.5,1.5) .. (-1,2);
		
\node at (0,-3.4){Type 2.03};			
		
	\end{tikzpicture}
\end{minipage}%
\begin{minipage}[t]{0.2\linewidth} 
\centering
\begin{tikzpicture}[xscale=0.6,yscale=0.6]
			
	\draw[blue,thick] (0,2)\nn to (0,2)\nn;

			\draw[blue,thick] (0,-2)\nn to (0,-2)\nn;
				\draw[blue,thick] (0,0)\nn to (0,0)\nn;

\draw[thick,color=blue] (0,-2) .. controls (-2,-1) and (-2,1) .. (0,2);	
\draw[thick,color=blue] (0,-2) .. controls (2,-1) and (2,1) .. (0,2);	

\draw[thick,color=blue] (0,-2) .. controls (-0.5,-1.5) and (-0.5,-0.5) .. (0,0);	
\draw[thick,color=blue] (0,-2) .. controls (0.5,-1.5) and (0.5,-0.5) .. (0,0);	
		
\draw[blue](0,-2)node[below]{$P_1$}(0,2)node[above]{$P_2$}(0,0)node[left]{$P$};	
		
\draw[blue](0.2,-0.3)node[rotate=20]{$\bowtie$};

\draw[thick,color=red] (1.4,-1.8) .. controls (0.2,-1.3) and (-0.2,-1.3) .. (-1.4,-1.8);	
		
\node at (0,-3.4){Type 2.04};	

\end{tikzpicture}
\end{minipage}%

\begin{minipage}[t]{0.2\linewidth} 
\centering
\begin{tikzpicture}[xscale=0.6,yscale=0.6]
			
		\draw[thick,red] (0,0)\nn to (0,2)\nn;
			\draw[blue,thick] (0,-2)\nn to (0,-2)\nn;
				\draw[blue,thick] (0,0)\nn to (0,0)\nn;
				
				\draw[blue,thick] (0,2)\nn to (0,2)\nn;
				
\draw[thick,color=blue] (0,-2) .. controls (-2,-1) and (-2,1) .. (0,2);	
\draw[thick,color=blue] (0,-2) .. controls (2,-1) and (2,1) .. (0,2);	

\draw[thick,color=blue] (0,-2) .. controls (-0.5,-1.5) and (-0.5,-0.5) .. (0,0);	
\draw[thick,color=blue] (0,-2) .. controls (0.5,-1.5) and (0.5,-0.5) .. (0,0);	
		
\draw[blue](0,-2)node[below]{$P_1$}(0,2)node[above]{$P_2$}(0,0)node[left]{$P$};

\draw[blue](0.2,-0.3)node[rotate=20]{$\bowtie$};
\node at (0,-3.4){Type 2.05};
\end{tikzpicture}

\end{minipage}
\begin{minipage}[t]{0.2\linewidth} 
\centering
\begin{tikzpicture}[xscale=0.6,yscale=0.6]

			\draw[thick,color=red] (0,2) .. controls (-0.8,1) and (-0.5,-0.3) .. (0,-0.3);						
\draw[thick,color=red] (0,2) .. controls (0.8,1) and (0.5,-0.3) .. (0,-0.3);	

			\draw[blue,thick] (0,-2)\nn to (0,-2)\nn;
				\draw[blue,thick] (0,0)\nn to (0,0)\nn;
				
				\draw[blue,thick] (0,2)\nn to (0,2)\nn;
				
\draw[thick,color=blue] (0,-2) .. controls (-2,-1) and (-2,1) .. (0,2);	
\draw[thick,color=blue] (0,-2) .. controls (2,-1) and (2,1) .. (0,2);	

\draw[thick,color=blue] (0,-2) .. controls (-0.5,-1.5) and (-0.5,-0.5) .. (0,0);	
\draw[thick,color=blue] (0,-2) .. controls (0.5,-1.5) and (0.5,-0.5) .. (0,0);	
		
\draw[blue](0,-2)node[below]{$P_1$}(0,2)node[above]{$P_2$}(0,0)node[left]{$P$};

\draw[blue](0.2,-0.3)node[rotate=20]{$\bowtie$};

\node at (0,-3.4){Type 2.06};

\end{tikzpicture}
\end{minipage}
\begin{minipage}[t]{0.2\linewidth}
\begin{tikzpicture}[xscale=0.6,yscale=0.6]
		\draw[blue,thick] (0,2)\nn to (0,2)\nn;

			\draw[blue,thick] (0,-2)\nn to (0,-2)\nn;
				\draw[blue,thick] (0,0)\nn to (0,0)\nn;

\draw[thick,color=blue] (0,-2) .. controls (-2,-1) and (-2,1) .. (0,2);	
\draw[thick,color=blue] (0,-2) .. controls (2,-1) and (2,1) .. (0,2);	

\draw[thick,color=blue] (0,-2) .. controls (-0.5,-1.5) and (-0.5,-0.5) .. (0,0);	
\draw[thick,color=blue] (0,-2) .. controls (0.5,-1.5) and (0.5,-0.5) .. (0,0);	
		
\draw[blue](0,-2)node[below]{$P_1$}(0,2)node[above]{$P_2$}(0,0)node[left]{$P$};	
\draw[blue](0.2,-0.3)node[rotate=20]{$\bowtie$};

\draw[thick,color=red] (1.8,-0.5) .. controls (-1,-2) and (-2,-1) .. (0,2);
\node at (0,-3.4){Type 2.07};
		\end{tikzpicture}

\end{minipage}
\begin{minipage}[t]{0.2\linewidth} 
\centering
\begin{tikzpicture}[xscale=0.6,yscale=0.6]
		
			\draw[blue,thick] (0,-2)\nn to (0,-2)\nn;
				\draw[blue,thick] (0,0)\nn to (0,0)\nn;
				
				\draw[blue,thick] (0,2)\nn to (0,2)\nn;
				
\draw[thick,color=blue] (0,-2) .. controls (-2,-1) and (-2,1) .. (0,2);	
\draw[thick,color=blue] (0,-2) .. controls (2,-1) and (2,1) .. (0,2);	

\draw[thick,color=blue] (0,-2) .. controls (-0.5,-1.5) and (-0.5,-0.5) .. (0,0);	
\draw[thick,color=blue] (0,-2) .. controls (0.5,-1.5) and (0.5,-0.5) .. (0,0);	
		
\draw[blue](0,-2)node[below]{$P_1$}(0,2)node[above]{$P_2$}(0,0)node[left]{$P$};	
		
\draw[blue](0.2,-0.3)node[rotate=20]{$\bowtie$};


 \draw[thick,color=red] (-1.8,-0.5) .. controls (1,-2) and (1.8,-1) .. (0,2);

\node at (0,-3.4){Type 2.08};
\end{tikzpicture}
\end{minipage}

\begin{minipage}[t]{0.2\linewidth} 
\centering
\begin{tikzpicture}[xscale=0.6,yscale=0.6]
		
			\draw[blue,thick] (0,-2)\nn to (0,-2)\nn;
				\draw[blue,thick] (0,0)\nn to (0,0)\nn;
				
		\draw[blue,thick] (0,2)\nn to (0,2)\nn;	
				
\draw[thick,color=blue] (0,-2) .. controls (-2,-1) and (-2,1) .. (0,2);	
\draw[thick,color=blue] (0,-2) .. controls (2,-1) and (2,1) .. (0,2);	

\draw[thick,color=blue] (0,-2) .. controls (-0.5,-1.5) and (-0.5,-0.5) .. (0,0);	
\draw[thick,color=blue] (0,-2) .. controls (0.5,-1.5) and (0.5,-0.5) .. (0,0);	
		
\draw[blue](0,-2)node[below]{$P_1$}(0,2)node[above]{$P_2$}(0,0)node[left]{$P$};	

\draw[blue](0.2,-0.3)node[rotate=20]{$\bowtie$};	

\draw[thick,color=red] (-2,0) .. controls (1,-2) and (2,-1) .. (-1.5,2);

\node at (0,-3.4){Type 2.09};
	\end{tikzpicture}
	
\end{minipage}%
\begin{minipage}[t]{0.2\linewidth} 
\centering
\begin{tikzpicture}[xscale=0.6,yscale=0.6]
			\draw[blue,thick] (0,2)\nn to (0,2)\nn;

			\draw[blue,thick] (0,-2)\nn to (0,-2)\nn;
				\draw[blue,thick] (0,0)\nn to (0,0)\nn;

\draw[thick,color=blue] (0,-2) .. controls (-2,-1) and (-2,1) .. (0,2);	
\draw[thick,color=blue] (0,-2) .. controls (2,-1) and (2,1) .. (0,2);	

\draw[thick,color=blue] (0,-2) .. controls (-0.5,-1.5) and (-0.5,-0.5) .. (0,0);	
\draw[thick,color=blue] (0,-2) .. controls (0.5,-1.5) and (0.5,-0.5) .. (0,0);	
		
\draw[blue](0,-2)node[below]{$P_1$}(0,2)node[above]{$P_2$}(0,0)node[left]{$P$};	
\draw[blue](0.2,-0.3)node[rotate=20]{$\bowtie$};	

\draw[thick,color=red] (2,0) .. controls (-1,-2) and (-2,-1) .. (1.5,2);

\node at (0,-3.4){Type 2.10};
		\end{tikzpicture}
\end{minipage}
 \begin{minipage}[t]{0.2\linewidth} 
\centering
\begin{tikzpicture}[xscale=0.6,yscale=0.6]
			
	\draw[blue,thick] (0,2)\nn to (0,2)\nn;

			\draw[blue,thick] (0,-2)\nn to (0,-2)\nn;
				\draw[blue,thick] (0,0)\nn to (0,0)\nn;

\draw[thick,color=blue] (0,-2) .. controls (-2,-1) and (-2,1) .. (0,2);	
\draw[thick,color=blue] (0,-2) .. controls (2,-1) and (2,1) .. (0,2);	

\draw[thick,color=blue] (0,-2) .. controls (-0.5,-1.5) and (-0.5,-0.5) .. (0,0);	
\draw[thick,color=blue] (0,-2) .. controls (0.5,-1.5) and (0.5,-0.5) .. (0,0);	
		
\draw[blue](0,-2)node[below]{$P_1$}(0,2)node[above]{$P_2$}(0,0)node[left]{$P$};	
\draw[blue](0.2,-0.3)node[rotate=20]{$\bowtie$};	

	 \draw[thick,red] (0,0) to (-2,1.5);

\node at (0,-3.4){Type 2.11};

	\end{tikzpicture}
\end{minipage}%
\begin{minipage}[t]{0.2\linewidth} 
\centering
\begin{tikzpicture}[xscale=0.6,yscale=0.6]
			
	\draw[blue,thick] (0,2)\nn to (0,2)\nn;
			
			 \draw[thick,red] (0,0) to (2,1.5);	
			
			\draw[blue,thick] (0,-2)\nn to (0,-2)\nn;
				\draw[blue,thick] (0,0)\nn to (0,0)\nn;

\draw[thick,color=blue] (0,-2) .. controls (-2,-1) and (-2,1) .. (0,2);	
\draw[thick,color=blue] (0,-2) .. controls (2,-1) and (2,1) .. (0,2);	

\draw[thick,color=blue] (0,-2) .. controls (-0.5,-1.5) and (-0.5,-0.5) .. (0,0);	
\draw[thick,color=blue] (0,-2) .. controls (0.5,-1.5) and (0.5,-0.5) .. (0,0);	
		
\draw[blue](0,-2)node[below]{$P_1$}(0,2)node[above]{$P_2$}(0,0)node[left]{$P$};	
			
\draw[blue](0.2,-0.3)node[rotate=20]{$\bowtie$};	

\node at (0,-3.4){Type 2.12};

	\end{tikzpicture}
\end{minipage}%

  \caption{Irreducible arc segments in a type II tile}\label{fig:irreducible-arc-segments-type-II}
  \end{figure}

Fix  a tagged triangulation $\fT$ of $(\fS,\fM)$. Every permissible  tagged arc  is decomposed by $\fT$ into segments. The equivalence of arcs induces an equivalence relation on  segments lying in the same region of $\fT$. An equivalence class of a segment is called an {\em irreducible arc segment} (with respect to $\fT$). All possible irreducible are segments in a type I tile are listed in Figure \ref{fig:type-I-irred-arc-seg}, while all possible irreducible arc segments in a type II tile are listed in Figure \ref{fig:irreducible-arc-segments-type-II}, by omitting the possible tags. More precisely, only for type 1.1, 2.05, 2.11 and 2.12, the end incident to $P$ may be tagged notched, and in type 1.3, $\eta$ has a different tag with $\ba$ at the end incident to $P$.

Let  $\cm$ be a finite permissible multiset with respect to $\mathbf{T}$ and $\Delta$ a tile. We denote by
$\Arcseg_{\Delta}(\cM)$ the finite multiset of all  the irreducible arc segments  of arcs in $\cM$ lying in $\Delta$.

\section{Intersection vector determines irreducible arc segments}\label{s:local-arc-segmet}
Recall that $(\fS,\fM)$ is a disc with at most one puncture, which we denote by $P$ if it exists. Let $\fT$ be a tagged triangulation of $(\fS,\fM)$. The aim of this section is to show that the intersection vector of a finite permissible multiset $\cm$ of tagged arcs determines the irreducible arc segments of $\cm$. Namely,
\begin{thm}\label{thm:intersection-arc-segments}
    Let $\cm$ and $\cn$ be two finite permissible multisets of tagged arcs with respect to $\fT$ such that $\Intv_\fT(\cm)=\Intv_{\fT}(\cn)$. Then we have $\Arcseg_{\Delta}(\cm)=\Arcseg_{\Delta}(\cn)$ for each tile $\Delta$ of $(\fS,\fM)$.
\end{thm}
Theorem \ref{thm:intersection-arc-segments} is a direct consequence of Lemma \ref{l:tri-gons}, Proposition \ref{prop:traingle-puncture} and Proposition \ref{prop:1-puncture}. 
 Throughout this section, we always assume that $\cm$ and $\cn$ satisfy the condition of Theorem \ref{thm:intersection-arc-segments}.
For an arc $\gamma$ and a positive integer $m$, we denote by $\gamma^{(m)}$ the multiset that contains $m$ copies of $\gamma$. We also denote by $\gamma^{(0)}=\emptyset$.


\subsection{Triangle with all its vertices on the boundary}\label{ss:proper-triangle}
Let $\Delta$ be a triangle as in Figure \ref{f:basic triangle} and $\cc$ a finite permissible multiset of tagged arcs with respect to $\fT$. We denote by
\begin{itemize}
    \item $\Arc_\Delta(\cc, \angle \eta_i)$ the multiset of irreducible arc segments in $\Arcs_\Delta(\cc)$ that cross $\angle \eta_i$, for $i=1,2,3$.
    \item $\Arc_\Delta(\cc,P_i,\ba_i)$ the multiset of irreducible arc segments in $\Arc_\Delta(\cc)$ that have one endpoint at $P_i$ and the other in interior of $\ba_i$.
    \item $l_{\eta_i}^\cc=|\Arc_\Delta(\cc, \angle \eta_i)|$ the cardinality of $\Arc_\Delta(\cc, \angle \eta_i)$.
    \item $l_{P_i\ba_i}^\cc=|\Arc_\Delta(\cc,P_i,\ba_i)|$ the cardinality of $\Arc_\Delta(\cc,P_i,\ba_i)$.
\end{itemize}
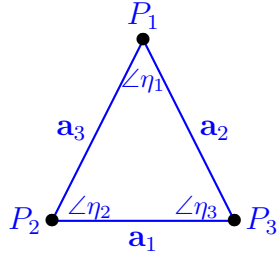
\begin{figure}[ht]

\begin{tikzpicture}[xscale=0.6,yscale=0.6]

		\draw[blue,thick] (-2,-2)\nn to (2,-2)\nn;
			
		\draw[blue,thick] (-2,-2)\nn to (0,2)\nn;
			\draw[blue,thick] (2,-2)\nn to (0,2)\nn;

		\draw[blue](-2,-2)node[left]{$P_2$}(-1.9,-1.7)node[right]{\footnotesize$\angle\eta_2$}(1.9,-1.7)node[left]{\footnotesize$\angle\eta_3$}(2,-2)node[right]{$P_3$}(0,2)node[above]{$P_1$}(0,1.6)node[below]{\footnotesize$\angle\eta_1$};	
		\draw[blue](0,-2)node[below]{$\ba_1$}(1,0)node[right]{$\ba_2$}(-1,0)node[left]{$\ba_3$};			
	
	\end{tikzpicture}
\caption{Triangle with all its vertices on the boundary}\label{f:basic triangle}
\end{figure}
\begin{lem}\label{l:tri-gons}
Let $\Delta$ be  a triangle with all its vertices on the boundary, then \[\Arcseg_{\Delta}(\cm)=\Arcseg_{\Delta}(\cn).\]
\end{lem}

\begin{proof}

Let $\Delta$ be a triangle of type I as described in Figure~\ref{f:basic triangle} such that $P_1,P_2,P_3$ are  boundary marked points. It suffices to show that $l_{\eta_i}^\cm=l_{\eta_i}^\cn$ and $l_{P_i\ba_i}^\cm=l_{P_i\ba_i}^\cn$ for $i=1,2,3$.

Recall that $\fT$ contains no loops, $\cm\cap \fT=\emptyset=\cn\cap\TT$ and
 $P_1,P_2,P_3$ are  on the boundary. It follows  that \[\Int^A(\ba_i|\cm)=\Int(\ba_i|\cm)=\Int(\ba_i|\cn)=\Int^A(\ba_i|\cn)\] for each $1\leq i\leq 3$.
On the other hand, we have
\[
\begin{cases}
\Int^A(\ba_1|\cm)=l_{\eta_2}^\cm+l_{\eta_3}^\cm+l_{P_{1}\ba_1}^\cm\\
\Int^A(\ba_2|\cm)=l_{\eta_1}^\cm+l_{\eta_3}^\cm+l_{P_{2}\ba_2}^\cm\\
\Int^A(\ba_3|\cm)=l_{\eta_1}^\cm+l_{\eta_2}^\cm+l_{P_{3}\ba_3}^\cm
\end{cases}\quad \text{and}\quad \begin{cases}
\Int^A(\ba_1|\cn)=l_{\eta_2}^\cn+l_{\eta_3}^\cn+l_{P_{1}\ba_1}^\cn\\
\Int^A(\ba_2|\cn)=l_{\eta_1}^\cn+l_{\eta_3}^\cn+l_{P_{2}\ba_2}^\cn\\
\Int^A(\ba_3|\cn)=l_{\eta_1}^\cn+l_{\eta_2}^\cn+l_{P_{3}\ba_3}^\cn
\end{cases}
\]
by definition. It follows that 
\begin{eqnarray*}
\Int^A(\ba_2|\cm)+\Int^A(\ba_3|\cm)-\Int^A(\ba_1|\cm)&=&2l_{\eta_1}^\cm+l_{P_2\ba_2}^\cm+l_{P_3\ba_3}^\cm-l_{P_1\ba_1}^\cm,\\
\Int^A(\ba_2|\cn)+\Int^A(\ba_3|\cn)-\Int^A(\ba_1|\cn)&=&2l_{\eta_1}^\cn+l_{P_2\ba_2}^\cn+l_{P_3\ba_3}^\cn-l_{P_1\ba_1}^\cn.
\end{eqnarray*}
Consequently,
\begin{eqnarray}\label{eq:triangle-intersection-equality}
2l_{\eta_1}^\cm+l_{P_2\ba_2}^\cm+l_{P_3\ba_3}^\cm-l_{P_1\ba_1}^\cm=2l_{\eta_1}^\cn+l_{P_2\ba_2}^\cn+l_{P_3\ba_3}^\cn-l_{P_1\ba_1}^\cn.
\end{eqnarray}
We first note that if $l_{P_1\ba_1}^\cm\neq 0$, then $l_{\eta_1}^\cm=l_{P_2\ba_2}^\cm=l_{P_3\ba_3}^\cm=0$ since any two tagged arcs in $\cm$ are compatible. Similarly, if  $l_{P_1\ba_1}^\cn\neq 0$, then  $l_{\eta_1}^\cn=l_{P_2\ba_2}^\cn=l_{P_3\ba_3}^\cn=0$. According to (\ref{eq:triangle-intersection-equality}), we conclude that $l_{P_1\ba_1}^\cm\neq 0$ if and only if $l_{P_1\ba_1}^\cn\neq 0$. By the symmetry of the triangle $\Delta$, we conclude that $l_{P_i\ba_i}^\cm\neq 0$ if and only if $l_{P_i\ba_i}^\cn\neq 0$ for $i=1,2,3$. We separate the remaining proof by discussing whether there exists an $i$ such that $l_{P_i\ba_i}^\cm\neq 0$.

\noindent{\bf Case 1}: $l_{P_i\ba_i}^\cm\neq 0$ for some $i$. Without loss of generality, we may assume $i=1$. It follows that $l_{P_1\ba_1}^\cn\neq 0$ and $l_{\eta_1}^\cm=l_{P_2\ba_2}^\cm=l_{P_3\ba_3}^\cm=l_{\eta_1}^\cn=l_{P_2\ba_2}^\cn=l_{P_3\ba_3}^\cn=0$ (cf. the first graph of Figure~\ref{f:Cases of no edge of delta is loop}). By definition, we have
\[
\begin{cases}
\Int^A(\ba_1|\cm)=l_{\eta_2}^\cm+l_{\eta_3}^\cm+l_{P_{1}\ba_1}^\cm\\
\Int^A(\ba_2|\cm)=l_{\eta_3}^\cm\\
\Int^A(\ba_3|\cm)=l_{\eta_2}^\cm
\end{cases}\quad \text{and} \quad   \begin{cases}
\Int^A(\ba_1|\cn)=l_{\eta_2}^\cn+l_{\eta_3}^\cn+l_{P_{1}\ba_1}^\cn\\
\Int^A(\ba_2|\cn)=l_{\eta_3}^\cn\\
\Int^A(\ba_3|\cn)=l_{\eta_2}^\cn
\end{cases}.
\]
Since $\Int^A(\ba_i|\cm)=\Int^A(\ba_i|\cn)$, we obtain 
\[
\left\{
\begin{array}{ccc}
l_{\eta_2}^\cm+l_{\eta_3}^\cm+l_{P_{1}\ba_1}^\cm&=&l_{\eta_2}^\cn+l_{\eta_3}^\cn+l_{P_{1}\ba_1}^\cn\\
                l_{\eta_3}^\cm&=&l_{\eta_3}^\cn\\
               l_{\eta_2}^\cm&=&l_{\eta_2}^\cn
\end{array}.
\right.
\]
Consequently,  $l_{\eta_i}^\cm=l_{\eta_i}^\cn$ and $l_{P_i\ba_i}^\cm=l_{P_i\ba_i}^\cn$ for any $1\leq i\leq 3$.



\noindent{\bf Case 2:}  $l_{P_i\ba_i}^\cm= 0$ for all $1\leq i\leq 3$. Hence $l_{P_i\ba_i}^\cn= 0$ for all $1\leq i\leq 3$.
 It follows that
\[
\begin{cases}
\Int^A(\ba_1|\cm)=l_{\eta_2}^\cm+l_{\eta_3}^\cm\\
\Int^A(\ba_2|\cm)=l_{\eta_1}^\cm+l_{\eta_3}^\cm\\
\Int^A(\ba_3|\cm)=l_{\eta_1}^\cm+l_{\eta_2}^\cm
\end{cases}\quad \text{and} \quad  \begin{cases}
\Int^A(\ba_1|\cn)=l_{\eta_2}^\cn+l_{\eta_3}^\cn\\
\Int^A(\ba_2|\cn)=l_{\eta_1}^\cn+l_{\eta_3}^\cn\\
\Int^A(\ba_3|\cn)=l_{\eta_1}^\cn+l_{\eta_2}^\cn
\end{cases}.
\] 
(cf. the second graph of Figure~\ref{f:Cases of no edge of delta is loop}). Again, by  $\Int^A(\ba_i|\cm)=\Int^A(\ba_i|\cn)$, we have
\[
\left\{
\begin{array}{ccc}
l_{\eta_2}^\cm+l_{\eta_3}^\cm&=&l_{\eta_2}^\cn+l_{\eta_3}^\cn\\
l_{\eta_1}^\cm+l_{\eta_3}^\cm&=&l_{\eta_1}^\cn+l_{\eta_3}^\cn\\
l_{\eta_1}^\cm+l_{\eta_2}^\cm&=&l_{\eta_1}^\cn+l_{\eta_2}^\cn
\end{array}.
\right.
\]
Consequently, $l_{\eta_i}^\cm=l_{\eta_i}^\cn$ and $l_{P_i\ba_i}^\cm=l_{P_i\ba_i}^\cn$ for any $1\leq i\leq 3$.

	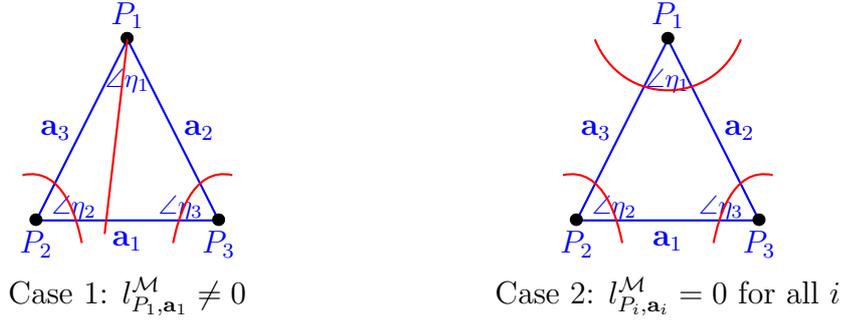
\begin{figure}[ht]

	\begin{minipage}{0.4\linewidth} 
		\begin{tikzpicture}[xscale=0.6,yscale=0.6]
			
			\draw[blue,thick] (-2,-2)\nn to (2,-2)\nn;
			
		\draw[blue,thick] (-2,-2)\nn to (0,2)\nn;
			\draw[blue,thick] (2,-2)\nn to (0,2)\nn;

		\draw[blue](-2,-2)node[below]{$P_2$}(2,-2)node[below]{$P_3$}(0,2)node[above]{$P_1$};	
		\draw[blue](-1.9,-1.7)node[right]{\footnotesize$\angle\eta_2$}(1.9,-1.7)node[left]{\footnotesize$\angle\eta_3$}(0,1.6)node[below]{\footnotesize$\angle\eta_1$};

		\draw[blue](0,-2)node[below]{$\ba_1$}(1,0)node[right]{$\ba_2$}(-1,0)node[left]{$\ba_3$};			
		
	\draw[red,thick] (0,2) to (-0.5,-2.3);
\draw[thick,color=red] (-2.3,-1) .. controls (-1.4,-0.8) and (-1.1,-2) .. (-1,-2.5);	
\draw[thick,color=red] (2.3,-1) .. controls (1.4,-0.8) and (1.1,-2) .. (1,-2.5);	
\draw[thick](0,-3)node[below]{Case 1: $l_{P_1,\ba_1}^\cm\neq 0$};

		\end{tikzpicture}
		\end{minipage}
\begin{minipage}{0.4\linewidth}
		\begin{tikzpicture}[xscale=0.6,yscale=0.6]
			
			\draw[blue,thick] (-2,-2)\nn to (2,-2)\nn;
			
		\draw[blue,thick] (-2,-2)\nn to (0,2)\nn;
			\draw[blue,thick] (2,-2)\nn to (0,2)\nn;

		\draw[blue](-2,-2)node[below]{$P_2$}(2,-2)node[below]{$P_3$}(0,2)node[above]{$P_1$};	
		\draw[blue](-1.9,-1.7)node[right]{\footnotesize$\angle\eta_2$}(1.9,-1.7)node[left]{\footnotesize$\angle\eta_3$}(0,1.6)node[below]{\footnotesize$\angle\eta_1$};

		\draw[blue](0,-2)node[below]{$\ba_1$}(1,0)node[right]{$\ba_2$}(-1,0)node[left]{$\ba_3$};

\draw[thick,color=red] (-2.3,-1) .. controls (-1.4,-0.8) and (-1.1,-2) .. (-1,-2.5);	

\draw[thick,color=red] (2.3,-1) .. controls (1.4,-0.8) and (1.1,-2) .. (1,-2.5);	
	
\draw[thick,color=red] (-1.6,2) .. controls (-1,0.5) and (1,0.5) .. (1.6,2);	
\draw[thick](0,-3)node[below]{Case 2: $l_{P_i,\ba_i}^\cm=0$ for all $i$};				
		\end{tikzpicture}
		
		\end{minipage}
	
	\caption{$\Delta$ is a triangle with all its vertices on the boundary}\label{f:Cases of no edge of delta is loop}
	\end{figure}

\end{proof}

\subsection{Triangle with the puncture $P$ as a vertex}
Assume that $\fT$ has a triangle with the puncture $P$ as a vertex. Then all the triangles with $P$ as a common vertex are glued together as a polygon (cf. Figure \ref{f:Delta_i}). Let $\Delta_1,\dots, \Delta_s$ be all the triangles with $P$ as a vertex. We may assume that $\Delta_i$ and $\Delta_{i+1}$ share a common edge $\ba_{i+1}$ with endpoints $P$ and $P_{i+1}$ for $1\leq i\leq s-1$, $\Delta_s$ shares an edge $\ba_1$ with $\Delta_1$ whose endpoints are $P$ and $P_1$. 
For each $\Delta_i$($1\leq i\leq s$), we denote 
\begin{itemize}
    \item $b_i$: the remaining edge of $\Delta_i$;
    \item $\angle \eta_i$: the angle  with vertex $P$;
    \item $\angle \delta_i$: the angle with vertex $P_i$;
    \item $\angle \sigma_i$:  the angle  with vertex $P_{i+1}$ (cf. Figure~\ref{f:Delta_i}).
\end{itemize}
Note that all the arcs $\ba_1,\dots,\ba_s$ have the same tags at the ends incident to $P$. On the other hand, since $\cm$ and $\cn$ are permissible, all the arcs in $\cm$ (resp. $\cn$) incident to $P$ also have the same tags at the ends incident to $P$. However, changing the tags of $\ba_1,\dots, \ba_s$ and all the arcs in $\cm$ (resp. $\cn$)  at the ends incident to $P$ do not change the intersection number. Hence we may assume that $\ba_1,\dots, \ba_s$ are tagged plain at the ends incident to $P$.

\begin{figure}
    \centering
    
    \begin{tikzpicture}[xscale=0.6,yscale=0.6]
			
			\draw[blue,thick] (-2,-2)\nn to (2,-2)\nn;
			
		\draw[blue,thick] (-2,-2)\nn to (0,2)\nn;
			\draw[blue,thick] (2,-2)\nn to (0,2)\nn;
			\draw[blue,thick] (0,2)\nn to (4,2)\nn;
			\draw[blue,thick] (0,2)\nn to (-4,2)\nn;
			\draw[blue,thick] (-2,-2)\nn to (-4,2)\nn;
			\draw[blue,thick] (2,-2)\nn to (4,2)\nn;
		
		\draw[blue](-2,-2)node[left]{$P_i$}(-1.9,-1.7)node[right]{\footnotesize$\angle\delta_i$}(1.9,-1.7)node[left]{\footnotesize$\angle \sigma_i$}(2,-2)node[right]{$P_{i+1}$}(0,2)node[above]{$P$}(0,1.6)node[below]{\footnotesize$\angle\eta_i$};	
		\draw[blue](0,-2)node[below]{$b_i$}(1,0)node[right]{$\ba_{i+1}$}(-1,0)node[left]{$\ba_i$};			
			
		\draw[blue](-3,0)node[left]{$b_{i-1}$}(-2,2)node[above]{$\ba_{i-1}$}(2,2)node[above]{$\ba_{i+2}$}(3,0)node[right]{$b_{i+1}$};		
		\draw[red](0,0)node[below]{$\Delta_i$}(-2,0.5)node[above]{$\Delta_{i-1}$}(2,0.3)node[above]{$\Delta_{i+1}$};

		\draw[thick,color=blue,dotted] (-2,3) .. controls (-1,4) and (1,4) .. (2,3);	
			
		\end{tikzpicture}

    \caption{$\Delta_i$}\label{f:Delta_i}
\end{figure}

Keep the notation as in Section \ref{ss:proper-triangle}. Denote by 
\begin{itemize}
    \item $l^\cm$ (resp. $l^\cn$) the cardinality of the multiset consisting of irreducible arc segments of arcs in $\cm$ (resp. $\cn$) which are tagged notched at $P$.
    \item $l_i^\cm$ (resp. $l_i^\cn$) the cardinality of the multiset consisting of  arcs in $\cm$ (resp. $\cn$) whose untagged version is $\ba_i$ but are notched at $P$.
\end{itemize} 
 Then for $\ast\in\{\cm,\cn\}$, we have
\begin{equation}\label{eq:intersection-number-ab}
   \begin{cases}
\Int(\ba_i|\ast)=l^{\ast}-l_i^{\ast}+l^{\ast}_{\eta_i}+l^{\ast}_{P_{i+1}\ba_i}+l^{\ast}_{\delta_i}\\
\Int(\ba_{i+1}|\ast)=l^{\ast}-l_{i+1}^{\ast}+l^{\ast}_{\eta_i}+l^{\ast}_{P_{i}\ba_{i+1}}+l^{\ast}_{\sigma_i}\\
\Int(b_i|\ast)=l^{\ast}_{Pb_i}+l^{\ast}_{\delta_i}+l^{\ast}_{\sigma_i}.
\end{cases} 
\end{equation}
Since $\ast$  is permissible, if $l^\ast=0$, then $l_i^\ast=0$ for $1\leq i\leq s$.
If $l^\ast>0$, then $l^\ast=\sum_{i=1}^s(l^\ast_{Pb_i}+l_i^\ast)$. In particular, 
$l^\ast_{Pb_i}+l_i^\ast\leq l^\ast$ for all $1\leq i\leq s$ if $l^\ast>0$.
\begin{remark}
  Let $\cc$ be a finite permissible multiset.  We remark that if $l^{\cc}=0$, there may exist $1\leq i\leq s$ such that $l^\cc_{Pb_i}>0$. In other words, there are irreducible arc segments of type $1.1$ in $\Delta_i$ which have plain tag at the end incident to $P$.
\end{remark}
\begin{lem}\label{lem:lm=ln=0}
    $l^\cm=0$ if and only if $l^\cn=0$.
\end{lem}
\begin{proof}
    It suffices to show that $l^\cm=0$ implies $l^\cn=0$. Assume that $l^\cm=0$ but $l^\cn>0$.  

For any $1\leq i\leq s$, by (\ref{eq:intersection-number-ab}), we have
    \begin{eqnarray}\label{eq:puncture-triangle-l-n>0}
        &&\Int(\ba_i|\cn)+\Int(\ba_{i+1}|\cn)-\Int(b_i|\cn)\notag\\
        &=&2l^\cn-l_i^\cn-l_{i+1}^\cn+2l_{\eta_i}^\cn+l_{P_{i+1}\ba_i}^\cn+l_{P_i\ba_{i+1}}^\cn-l_{Pb_i}^\cn\notag\\
       &=&(l^\cn-l_i^\cn-l_{Pb_i}^\cn)+(l^\cn-l_{i+1}^\cn)+2l_{\eta_i}^\cn+l_{P_{i+1}\ba_i}^\cn+l_{P_i\ba_{i+1}}^\cn\\
       &\geq&0.\notag
    \end{eqnarray}
Here, we use the convention that $\ba_{s+1}:=\ba_1$.
If there is an $i$ such that $l_{\eta_i}^\cm=l_{P_{i+1}\ba_i}^\cm=l_{P_i\ba_{i+1}}^\cm=0$, then by (\ref{eq:intersection-number-ab}) we have   \begin{align}\label{eq:puncture-triangle-l-m=0}
        \Int(\ba_i|\cm)+\Int(\ba_{i+1}|\cm)-\Int(b_i|\cm)=-l_{Pb_i}^\cm\leq 0.
    \end{align}
Since $\Int(\ba|\cm)=\Int(\ba|\cn)$ for any $\ba\in \fT$, by $(\ref{eq:puncture-triangle-l-n>0})$ and $(\ref{eq:puncture-triangle-l-m=0})$,
 we obtain
\[
(l^\cn-l_i^\cn-l_{Pb_i}^\cn)+(l^\cn-l_{i+1}^\cn)+2l_{\eta_i}^\cn+l_{P_{i+1}\ba_i}^\cn+l_{P_i\ba_{i+1}}^\cn=0,
\]
which contradicts $l^\cn>0$. 

Since $\cm$ consists of pairwise compatible tagged arcs,  there exists a $j$ such that $l_{\eta_j}^\cm=0$. Otherwise, there is a loop without endpoints enclosing $P$. If $l_{P_j\ba_{j+1}}^\cm=l_{P_{j+1}\ba_j}^\cm=0$, we obtain a contradiction by the above discussion. Hence, without loss of generality, we may assume that $l_{P_j\ba_{j+1}}^\cm\neq 0$. It follows that $l_{P_{j+1}\ba_j}^\cm=l_{\delta_j}^\cm=0$. Putting all of these into (\ref{eq:intersection-number-ab}), we obtain
 \[0=\Int(\ba_j|\cm)=\Int(\ba_j|\cn)=(l^\cn-l_j^\cn)+l_{\eta_j}^\cn+l_{P_{j+1}\ba_j}^\cn+l_{\delta_j}^\cn.\] 
Consequently, $l^\cn=l_j^\cn$ and $l_{\eta_j}^\cn=l_{\delta_j}^\cn=l_{P_{j+1}\ba_j}^\cn=0$.
Since $l_{\eta_j}^\cm=0$ and $l_{P_j\ba_{j+1}}^\cm\neq 0$, we conclude that $l_{\eta_{j-1}}^\cm=0$ and $l_{P_{j-1}\ba_j}^\cm=0$. If $l_{P_{j}\ba_{j-1}}^\cm=0$, we get a contradiction by the above discussion. Hence  $l_{P_{j}\ba_{j-1}}^\cm\neq 0$. If $s=2$,  there is  a once-puncture loop in $\cm$ which closely wraps around $\ba_j$, which contradicts the fact that $\cm$ is permissible. So we must have $s>2$. Then if for all $i\in \{1,\dots, s\}\backslash\{j-1,j\}$, $l_{\eta_i}^\cm\neq 0$,  there is  a once-puncture loop in $\cm$ which closely wraps around $\ba_j$, which again contradicts the fact that $\cm$ is permissible. Hence there is another $k\in \{1,\dots,s\}\backslash\{j-1,j\}$ such that $l_{\eta_k}^\cm=0$. Similar to the tile $\Delta_j$, we have either $l^\cn=l_k^\cn$ or $l^\cn=l_{k+1}^\cn$. Combining this with $l^\cn=l_j^\cn$, we obtain a contradiction with  $l^\cn=\sum_{i=1}^s(l_{Pb_i}^\cn+l_i^\cn)$. Hence $l^\cn=0$. This completes the proof.
\end{proof}

\begin{lem}\label{lem:lmln>0}
    Assume that $l^\cm>0$ and $l^\cn>0$, then $l^\cm=l^\cn$, $l_i^\cm=l_i^\cn$ and $l^\cm_{Pb_i}=l^\cn_{Pb_i}$ for each $1\leq i\leq s$.
\end{lem}
\begin{proof}

   Without loss of generality, we may assume that $l^\cm\geq l^\cn$. 
   We first claim that
   \begin{itemize}
       \item[(1)]  $l_i^\cm\geq l_i^\cn$ for each $1\leq i\leq s$.
       \item[(2)]  $l^\cm-l_i^\cm\geq l^\cn-l_i^\cn$ for each $1\leq i\leq s$.
       \item[(3)] $l^\cm_{Pb_i}\geq l^\cn_{Pb_i}$ for each $1\leq i\leq s$.
   \end{itemize}
   Suppose that there exists an $1\leq i\leq s$ such that $l_i^\cm<l_i^\cn$. In particular, $l_i^\cn>0$. Since $\cn$ is permissible, it follows that $l_{\eta_i}^\cn=l_{\delta_i}^\cn=l_{P_{i+1}\ba_i}^\cn=0$. Hence
  \[
       l^\cn-l_i^\cn=\Int(\ba_i|\cn)=\Int(\ba_i|\cm)
       =l^\cm-l_i^\cm+l_{\eta_i}^\cm+l_{\delta_i}^\cm+l_{P_{i+1}\ba_i}^\cm.
   \]
   Consequently, $l^\cn=l^\cm+(l_{\eta_i}^\cm+l_{\delta_i}^\cm+l_{P_{i+1}\ba_i}^\cm)+(l_i^\cn-l_i^\cm)>l^\cm$, a contradiction. Hence $l_i^\cn\leq l_i^\cm$ for all $1\leq i\leq s$.

  For $(2)$,  it is obvious if $l_i^\cm=0$. Assume that $l_i^\cm>0$. Then $\Int(\ba_i|\cm)=l^\cm-l_i^\cm$. On the other hand, 
  \[
  \Int(\ba_i|\cn)=l^\cn-l_i^\cn+l_{\eta_i}^\cn+l_{\delta_i}^\cn+l_{P_{i+1}\ba_i}^\cn.
  \]
  It follows that $(l^\cm-l_i^\cm)-(l^\cn-l_i^\cn)=l_{\eta_i}^\cn+l_{\delta_i}^\cn+l_{P_{i+1}\ba_i}^\cn\geq 0$. This proves the claim $(2)$.

   Now we prove the claim $(3)$.
    It suffices to assume that $l^\cn_{Pb_i}>0$. We have
    \begin{eqnarray*}
        2l^\cn-l_i^\cn-l_{i+1}^\cn-l_{Pb_i}^\cn&=& \Int(\ba_i|\cn)+\Int(\ba_{i+1}|\cn)-\Int(b_i|\cn)\\
        &=&\Int(\ba_i|\cm)+\Int(\ba_{i+1}|\cm)-\Int(b_i|\cm)\\
        &=&2l^\cm-l_i^\cm-l_{i+1}^\cm-l_{Pb_i}^\cm+(2l_{\eta_i}^\cm+l_{P_i\ba_{i+1}}^\cm+l_{P_{i+1}\ba_i}^\cm).
    \end{eqnarray*}
    It follows that 
    \[
    l_{Pb_i}^\cm-l_{Pb_i}^\cn\geq (l^\cm-l_i^\cm)-(l^{\cn}-l_i^\cn)+(l^\cm-l_{i+1}^\cm)-(l^{\cn}-l_{i+1}^\cn)\geq 0.
    \]

    Denote by
    \[\begin{array}{cc}
    \mathcal{X}_\cm:=\{i~|~1\leq i\leq s,l_{Pb_i}^\cm>0\}, & \mathcal{X}_\cn:=\{i~|~1\leq i\leq s,l_{Pb_i}^\cn>0\},\\
    \mathcal{Y}_\cm:=\{i~|~1\leq i\leq s, l_i^\cm>0\}, &\mathcal{Y}_\cn:=\{i~|~1\leq i\leq s, l_i^\cn>0\}.
    \end{array}
    \]
    According to $(1)$ and $(3)$, we know that $\mathcal{X}_\cn\subseteq \mathcal{X}_\cm$ and $\mathcal{Y}_{\cn}\subseteq \mathcal{Y}_\cm$.

    Let us first assume that $\mathcal{Y}_\cn\neq \emptyset$. In particular, for any $k\in \mathcal{Y}_\cn$, $l_k^\cm\geq l_k^\cn>0$. It follows that 
    \begin{align}\label{eq:lm-lk}
        l^\cm-l_k^\cm=\Int(\ba_k|\cm)=\Int(\ba_k|\cn)=l^\cn-l_k^\cn.
    \end{align}
    Note that $l^\cm=\sum_{i\in \mathcal{X}_\cm}l_{Pb_i}^\cm+\sum_{i\in \mathcal{Y}_\cm}l_i^\cm$ and $l^\cn=\sum_{i\in \mathcal{X}_\cn}l_{Pb_i}^\cn+\sum_{i\in \mathcal{Y}_\cn}l_i^\cn$. It follows that
    \begin{equation*}
        \sum_{i\in \mathcal{X}_\cm}l_{Pb_i}^\cm+\sum_{k\neq i\in \mathcal{Y}_\cm}l_i^\cm=\sum_{i\in \mathcal{X}_\cn}l_{Pb_i}^\cn+\sum_{k\neq i\in \mathcal{Y}_\cn}l_i^\cn.
    \end{equation*}
  According to $(1)$ and $(3)$, we conclude that
  \begin{align}\label{eq:yn-non-emptyset}
      \begin{cases}
      \text{ $\mathcal{X}_\cm=\mathcal{X}_\cn$ and $\mathcal{Y}_\cm=\mathcal{Y}_\cn$;}\\
      \text{ $l_{Pb_i}^\cm=l_{Pb_i}^\cn$ for all $1\leq i\leq s$;}\\
      \text{ $l_i^\cm=l_i^\cn$ for all $1\leq i\leq s$ but $i\neq k$.}
  \end{cases}
  \end{align}
 If  $|\mathcal{Y}_\cn|\geq 2$, we obtain that $l_i^\cm=l_i^\cn$ for all $1\leq i\leq s$ by applying (\ref{eq:yn-non-emptyset}) to different elements of $\mathcal{Y}_\cn$. Consequently, $l^\cm=l^\cn$. 
 
 Assume that $|\mathcal{Y}_\cn|=1$, say $\mathcal{Y}_\cn=\{k\}$. If $\mathcal{X}_\cm=\mathcal{X}_\cn\neq \emptyset$, then
  there exists a $j\in \mathcal{X}_\cm=\mathcal{X}_\cn$ such that $l_{Pb_j}^\cm=l_{Pb_j}^\cn>0$, we have
 \begin{eqnarray*}
     2l^\cm-l_j^\cm-l_{j+1}^\cm-l_{Pb_j}^\cm&=&\Int(\ba_j|\cm)+\Int(\ba_{j+1}|\cm)-\Int(b_j|\cm)\\
     &=&\Int(\ba_j|\cn)+\Int(\ba_{j+1}|\cn)-\Int(b_j|\cn)\\
     &=&2l^\cn-l_j^\cn-l_{j+1}^\cn-l_{Pb_j}^\cn.
 \end{eqnarray*}
  By (\ref{eq:lm-lk}) and (\ref{eq:yn-non-emptyset}),  we conclude that $l^\cm=l^\cn$.
  Now assume that $\mathcal{X}_\cm=\mathcal{X}_\cn=\emptyset$. Then $l^\cm=l_k^\cm$ and $l^\cn=l_k^\cn$. Denote by $\gamma_k$ the tagged arc such that $(\ba_k,\gamma_k)$ is a pair of conjugate arcs.
  Let $\cm_1=\cm\backslash {\gamma_k}^{(l_k^\cn)}$ and $\cn_1=\cn\backslash {\gamma_k}^{(l_k^\cn)}$. It follows that $\cm_1$ and $\cn_1$ are finite permissible multisets such that $\Intv_\fT(\cm_1)=\Intv_{\fT}(\cn_1)$. Furthermore, $l^{\cm_1}=l^\cm-l_k^\cn$ and $l^{\cn_1}=0$. By Lemma \ref{lem:lm=ln=0}, we get that $l^{\cm_1}=0$ and hence $l^\cm=l^\cn$.
    
   Now assume that $\mathcal{Y}_\cn=\emptyset$.  In this case, $l^\cn=\sum_{i\in \mathcal{X}_\cn}l_{Pb_i}^\cn>0$. It follows that $\mathcal{X}_\cn\neq \emptyset$.
     For any $k\in \mathcal{X}_\cn$, we have
    \[
    2l^\cm-l_k^\cm-l_{k+1}^\cm-l_{Pb_k}^\cm=2l^\cn-l_k^\cn-l_{k+1}^\cn-l_{Pb_k}^\cn=2l^\cn-l_{Pb_k}^\cn.
    \]
    It follows that \[
    \sum_{j\neq k,k+1}l_j^\cm+\sum_{j=1}^sl_j^\cm+2\sum_{j\in \mathcal{X}_\cm\backslash\mathcal{X}_\cn}l_{Pb_j}^\cm+2\sum_{k\neq j\in \mathcal{X}_\cn}(l_{Pb_j}^\cm-l_{Pb_j}^\cn)+(l_{Pb_k}^\cm-l_{Pb_k}^\cn)=0.
    \]
    Therefore $l_j^\cm=l_j^\cn=0$ for any $1\leq j\leq s$, $\mathcal{X}_\cm=\mathcal{X}_\cn$ and $l_{Pb_j}^\cm=l_{Pb_j}^\cn$ for any $j\in \mathcal{X}_\cn$, which implies that $l^\cm=l^\cn$. 
    \end{proof}

\begin{prop}\label{prop:traingle-puncture}
    Let $\Delta$ be a triangle with $P$ as a vertex. Then $\Arcs_\Delta(\cm)=\Arcs_{\Delta}(\cn)$.
\end{prop}
\begin{proof}
    Let $\Delta=\Delta_i$ as in Figure \ref{f:Delta_i}.
    If both $l^\cm=l^\cn=0$, then every arc in $\cm$ and $\cn$ is tagged plain at the end incident to $P$. Hence, we are in the situation  described in Section \ref{ss:proper-triangle}, and the result follows from Lemma \ref{l:tri-gons}. Now suppose $l^\cm>0$, then $l^\cn>0$ by Lemma \ref{lem:lm=ln=0}. Furthermore, we have $l^\cm=l^\cn$, $l_i^\cm=l_i^\cm$ and $l_{Pb_i}^\cm=l_{Pb_i}^\cn$ by Lemma \ref{lem:lmln>0}. Now, we apply the same proof as in Lemma \ref{l:tri-gons} to deduce the result.
\end{proof}

\subsection{$1$-puncture piece}
\begin{figure}[htpb]
		\begin{tikzpicture}[xscale=0.8,yscale=0.8]

		\draw[blue,thick] (0,2)\nn to (0,2)\nn;
			\draw[blue,thick] (0,-2)\nn to (0,-2)\nn;
				\draw[blue,thick] (0,0)\nn to (0,0)\nn;
\draw[thick,color=blue] (0,-2) .. controls (-2,-1) and (-2,1) .. (0,2);	
\draw[thick,color=blue] (0,-2) .. controls (2,-1) and (2,1) .. (0,2);	

\draw[thick,color=blue] (0,-2) .. controls (-0.5,-1.5) and (-0.5,-0.5) .. (0,0);	
\draw[thick,color=blue] (0,-2) .. controls (0.5,-1.5) and (0.5,-0.5) .. (0,0);	
		
\draw[blue](0,-2)node[below]{$P_1$}(0,2)node[above]{$P_2$}(0,0)node[above]{$P$};	
\draw[blue](-1.5,0)node[left]{$\ba_1$}(1.5,0)node[right]{$\ba_2$}(-0.3,-0.5)node[left]{$\ba_3$}(0.3,-0.5)node[right]{$\ba_4$};			
\draw[blue](0.2,-0.3)node[rotate=20]{$\bowtie$};			
		\end{tikzpicture}
		\caption{1-puncture piece}\label{f:1-puncture piece}
	\end{figure}
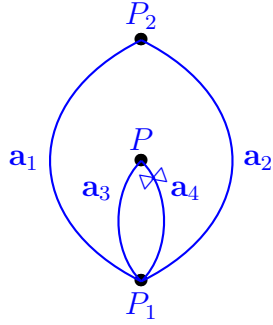

Let $\Delta$ be a $1$-puncture piece as described in Figure ~\ref{f:1-puncture piece}, where $P$ is the puncture and $P_1,P_2$ are different boundary marked points. 
Since  $\cm\cap \fT=\emptyset=\cn\cap \fT$ and $\ba_1,\dots,\ba_4$ are not loops, we obtain
\[
\Int(\ba_i|\cm)=\Int^A(\ba_i|\cm)+\Int^D(\ba_i|\cm)\quad \text{and}\quad \Int(\ba_i|\cn)=\Int^A(\ba_i|\cn)+\Int^D(\ba_i|\cn)
\]
for $1\leq i\leq 4$. It follows that 
\begin{eqnarray*}
    \Int^D(\ba_3|\cm)-\Int^D(\ba_4|\cm)&=&\Int(\ba_3|\cm)-\Int^A(\ba_3|\cm)-\Int(\ba_4|\cm)+\Int^A(\ba_4|\cm)\\
    &=&\Int(\ba_3|\cm)-\Int(\ba_4|\cm)\\
    &=&\Int(\ba_3|\cn)-\Int(\ba_4|\cn)\\
    &=&\Int(\ba_3|\cn)-\Int^A(\ba_3|\cn)-\Int(\ba_4|\cn)+\Int^A(\ba_4|\cn)\\
    &=&\Int^D(\ba_3|\cn)-\Int^D(\ba_4|\cn).
\end{eqnarray*}
Since $\cm$ and $\cn$ are permissible multisets, all arcs  in $\cm$ (resp. $\cn$) with $P$ as an endpoint are tagged in the same way at the end incident to $P$. In particular, $\Int^D(\ba_3|\cm)\times\Int^D(\ba_4|\cm)=0$ and $\Int^D(\ba_3|\cn)\times\Int^D(\ba_4|\cn)=0$. Consequently, all arcs with endpoints $P$ in $\cm$ are notched (resp. plain) at $P$ if and only if all arcs with endpoints $P$ in $\cn$ are notched (resp. plain) at $P$. Without loss of generality, we may assume that all arcs with endpoints $P$ in $\cm$ and $\cn$ are notched at $P$.

Denote by $l_{ij}^\cm$ (resp. $l_{ij}^\cn$) the cardinality of the sub-multiset of $\Arcs_\Delta(\cm)$ (resp. $\Arcs_\Delta(\cn)$) consisting of irreducible arc segments of type $2.ij$ (cf. Figure \ref{fig:irreducible-arc-segments-type-II}). Then $\Arcs_\Delta(\cm)$ and $\Arcs_\Delta(\cn)$ satisfy one of the $16$ cases in Table \ref{table:16-cases}.
\begin{table}[ht]
\centering
\begin{tabular}{|c|c|c|c|}
    \hline
    \multicolumn{2}{|c|}{(A): \quad $l_{01}^\ast+l_{02}^\ast+l_{03}^\ast\neq 0$ (cf. Figure \ref{f:Case A})}&\multicolumn{2}{|c|}{(B): \quad  $\begin{array}{c}
l_{01}^\ast+l_{02}^\ast+l_{03}^\ast=0\\
 l_{05}^\ast+l_{06}^\ast+l_{07}^\ast+l_{08}^\ast\neq 0
\end{array}$ (cf. Figure \ref{f:Case B})}\\ \hline
    (A1)& $l_{01}^\ast\neq 0$& (B1)& $l_{06}^\ast\neq 0$, $l_{07}^\ast\neq 0$\\ \hline
    (A2) & $l_{02}^\ast\neq 0$ &(B2) & $l_{06}^\ast\neq 0$, $l_{08}^\ast\neq 0$  \\ \hline
    (A3) & $l_{01}^\ast=l_{02}^\ast= 0$, $l_{03}^\ast\neq 0$, $l_{09}^\ast\neq 0$ &(B3) & $l_{06}^\ast\neq 0$, $l_{07}^\ast=l_{08}^\ast= 0$ \\ \hline
    (A4) & $l_{01}^\ast=l_{02}^\ast= 0$, $l_{03}^\ast\neq 0$, $l_{10}^\ast\neq 0$ & (B4) & $l_{06}^\ast= 0$, $l_{05}^\ast\neq 0$, $l_{07}^\ast\neq 0$ \\ \hline
    (A5) & $l_{01}^\ast=l_{02}^\ast=l_{09}^\ast=l_{10}^\ast= 0$, $l_{03}^\ast\neq 0$ & (B5) & $l_{06}^\ast= 0$, $l_{05}^\ast\neq 0$, $l_{08}^\ast\neq 0$  \\ \hline
    \multicolumn{2}{|c|}{(C): \quad $-l_{04}^\ast+\sum_{i=1}^8l_{0i}^\ast=0$(cf. Figure \ref{f:Case C})}&(B6) & $l_{06}^\ast=l_{07}^\ast=l_{08}^\ast= 0$, $l_{05}^\ast\neq 0$\\ \hline
     (C1)& $l_{09}^\ast\neq 0$&(B7) & $l_{05}^\ast=l_{06}^\ast= 0$, $l_{07}^\ast\neq 0$  \\ \hline
    (C2) & $l_{10}^\ast\neq 0$ &(B8) & $l_{05}^\ast=l_{06}^\ast= 0$, $l_{08}^\ast\neq 0$ \\ \hline
    (C3) & $l_{09}^\ast=l_{10}^\ast= 0$&& \\ \hline
\end{tabular}
\caption{Case (A)(B)(C), where $\ast\in\{\cm,\cn\}$}\label{table:16-cases}
\end{table}
Moreover, we have
\[
\begin{cases}
\Int(\ba_1|\ast)=l^{\ast}_{01}+l^{\ast}_{03}+l^{\ast}_{04}+l^{\ast}_{08}+2l^{\ast}_{09}+l^{\ast}_{11}\\
\Int(\ba_2|\ast)=l^{\ast}_{02}+l^{\ast}_{03}+l^{\ast}_{04}+l^{\ast}_{07}+2l^{\ast}_{10}+l^{\ast}_{12}\\ 
\Int(\ba_3|\ast)=l^{\ast}_{04}+l^{\ast}_{05}+l^{\ast}_{06}+l^{\ast}_{07}+l^{\ast}_{08}+l^{\ast}_{09}+l^{\ast}_{10}+l^{\ast}_{11}+l^{\ast}_{12}\\
\Int(\ba_4|\ast)=l^{\ast}_{04}+l^{\ast}_{06}+l^{\ast}_{07}+l^{\ast}_{08}+l^{\ast}_{09}+l^{\ast}_{10}.
\end{cases}
\]
for $\ast\in\{\cm,\cn\}$. Note that some irreducible arc segments of type $2.ij$ cannot appear together,  and all $l_{ij}^{\ast}$ that do not appear in Figure~\ref{f:Case A}, Figure~\ref{f:Case B} and Figure~\ref{f:Case C} are zeros.
\begin{figure}[ht]

\begin{minipage}[t]{0.3\linewidth} 
\centering
\begin{tikzpicture}[xscale=0.8,yscale=0.8]

			\draw[thick,red] (0,0) to (-2,0);	
			\draw[blue,thick] (0,-2)\nn to (0,-2)\nn;
				\draw[blue,thick] (0,0)\nn to (0,0)\nn;
				
		\draw[blue,thick] (0,2)\nn to (0,2)\nn;	
				
\draw[thick,color=blue] (0,-2) .. controls (-2,-1) and (-2,1) .. (0,2);	
\draw[thick,color=blue] (0,-2) .. controls (2,-1) and (2,1) .. (0,2);	

\draw[thick,color=blue] (0,-2) .. controls (-0.5,-1.5) and (-0.5,-0.5) .. (0,0);	
\draw[thick,color=blue] (0,-2) .. controls (0.5,-1.5) and (0.5,-0.5) .. (0,0);	
		
\draw[blue](0,-2)node[below]{$P_1$}(0,2)node[above]{$P_2$}(0,0)node[left]{$P$};	
\draw[blue](0.2,-0.3)node[rotate=20]{$\bowtie$};	
\draw[red](-0.4,0)node[rotate=90]{$\bowtie$};		

\draw[thick,color=red] (-1.6,-1.6) .. controls (1,-1) and (1,1) .. (-1.6,1.3);	
\draw[thick,color=red] (-1,2) .. controls (-0.5,1.5) and (0.5,1.5) .. (1,2);
\draw[thick,color=red] (0,-2) .. controls (1,-1) and (1.5,1) .. (-1.5,1.5);
\draw[blue](-0.9,0)node{$l_{11}^\ast$}(0,1.6)node{$l_{03}^\ast$}(0,0.6)node{$l_{09}^\ast$}(1,0)node{$l_{01}^\ast$};
\node at (0,-3.4){(A1): $l_{01}^\ast\neq 0$};
	\end{tikzpicture}
	
\end{minipage}%
\begin{minipage}[t]{0.3\linewidth} 
\centering
\begin{tikzpicture}[xscale=0.8,yscale=0.8]
			\draw[blue,thick] (0,2)\nn to (0,2)\nn;
			
			\draw[thick,red] (0,0) to (2,0);	
			\draw[blue,thick] (0,-2)\nn to (0,-2)\nn;
				\draw[blue,thick] (0,0)\nn to (0,0)\nn;

\draw[thick,color=blue] (0,-2) .. controls (-2,-1) and (-2,1) .. (0,2);	
\draw[thick,color=blue] (0,-2) .. controls (2,-1) and (2,1) .. (0,2);	

\draw[thick,color=blue] (0,-2) .. controls (-0.5,-1.5) and (-0.5,-0.5) .. (0,0);	
\draw[thick,color=blue] (0,-2) .. controls (0.5,-1.5) and (0.5,-0.5) .. (0,0);	
		
\draw[blue](0,-2)node[below]{$P_1$}(0,2)node[above]{$P_2$}(0,0)node[left]{$P$};	
\draw[blue](0.2,-0.3)node[rotate=20]{$\bowtie$};	
\draw[red](0.4,0)node[rotate=90]{$\bowtie$};		
\draw[thick,color=red] (1.6,-1.6) .. controls (-1,-1) and (-1,1) .. (1.6,1.3);	
\draw[thick,color=red] (1,2) .. controls (0.5,1.5) and (-0.5,1.5) .. (-1,2);
\draw[thick,color=red] (0,-2) .. controls (-1,-1) and (-1.5,1) .. (1.5,1.5);

\draw[blue](0.9,0)node{$l_{12}^\ast$}(0,1.6)node{$l_{03}^\ast$}(0.8,0.8)node{$l_{10}^\ast$}(-1,0)node{$l_{02}^\ast$};
\node at (0,-3.4){(A2): $ l_{02}^\ast\neq 0$};
		\end{tikzpicture}

\end{minipage}

\begin{minipage}[t]{0.3\linewidth} 
\centering
\begin{tikzpicture}[xscale=0.8,yscale=0.8]

			\draw[thick,red] (0,0) to (-2,1.5);	
			\draw[blue,thick] (0,-2)\nn to (0,-2)\nn;
				\draw[blue,thick] (0,0)\nn to (0,0)\nn;
				
		\draw[blue,thick] (0,2)\nn to (0,2)\nn;	
				
\draw[thick,color=blue] (0,-2) .. controls (-2,-1) and (-2,1) .. (0,2);	
\draw[thick,color=blue] (0,-2) .. controls (2,-1) and (2,1) .. (0,2);	

\draw[thick,color=blue] (0,-2) .. controls (-0.5,-1.5) and (-0.5,-0.5) .. (0,0);	
\draw[thick,color=blue] (0,-2) .. controls (0.5,-1.5) and (0.5,-0.5) .. (0,0);	
		
\draw[blue](0,-2)node[below]{$P_1$}(0,2)node[above]{$P_2$}(0,0)node[left]{$P$};	
\draw[blue](0.2,-0.3)node[rotate=20]{$\bowtie$};	
\draw[red](-0.4,0.3)node[rotate=45]{$\bowtie$};		

\draw[thick,color=red] (-1.4,-1.8) .. controls (-0.2,-1.3) and (0.2,-1.3) .. (1.4,-1.8);	
\draw[thick,color=red] (-1,2) .. controls (-0.5,1.5) and (0.5,1.5) .. (1,2);
\draw[thick,color=red] (-2,0) .. controls (1,-2) and (2,-1) .. (-1.5,2);
\draw[blue](-0.7,0.57)node{$l_{11}^\ast$}(0,1.6)node{$l_{03}^\ast$}(-1,-0.6)node{$l_{09}^\ast$}(0,-1.4)node{$l_{04}^\ast$};
\node at (0,-4){(A3):$\begin{array}{c}
 l_{01}^\ast=l_{02}^\ast=0\\
l_{03}^\ast\neq 0, l_{09}^\ast\neq 0
\end{array}$ };
	\end{tikzpicture}
	
\end{minipage}%
\begin{minipage}[t]{0.3\linewidth} 
\centering
\begin{tikzpicture}[xscale=0.8,yscale=0.8]
			\draw[blue,thick] (0,2)\nn to (0,2)\nn;
			
			\draw[thick,red] (0,0) to (2,1.5);	
			\draw[blue,thick] (0,-2)\nn to (0,-2)\nn;
				\draw[blue,thick] (0,0)\nn to (0,0)\nn;

\draw[thick,color=blue] (0,-2) .. controls (-2,-1) and (-2,1) .. (0,2);	
\draw[thick,color=blue] (0,-2) .. controls (2,-1) and (2,1) .. (0,2);	

\draw[thick,color=blue] (0,-2) .. controls (-0.5,-1.5) and (-0.5,-0.5) .. (0,0);	
\draw[thick,color=blue] (0,-2) .. controls (0.5,-1.5) and (0.5,-0.5) .. (0,0);	
		
\draw[blue](0,-2)node[below]{$P_1$}(0,2)node[above]{$P_2$}(0,0)node[left]{$P$};	
\draw[blue](0.2,-0.3)node[rotate=20]{$\bowtie$};	
\draw[red](0.4,0.3)node[rotate=-45]{$\bowtie$};		

\draw[thick,color=red] (1.4,-1.8) .. controls (0.2,-1.3) and (-0.2,-1.3) .. (-1.4,-1.8);	
\draw[thick,color=red] (1,2) .. controls (0.5,1.5) and (-0.5,1.5) .. (-1,2);
\draw[thick,color=red] (2,0) .. controls (-1,-2) and (-2,-1) .. (1.5,2);

\draw[blue](0.7,0.57)node{$l_{12}^\ast$}(0,1.6)node{$l_{03}^\ast$}(1,-0.6)node{$l_{10}^\ast$}(0,-1.4)node{$l_{04}^\ast$};
\node at (0,-4){(A4):
$\begin{array}{c}
 l_{01}^\ast=l_{02}^\ast=0\\
l_{03}^\ast\neq 0, l_{10}^\ast\neq 0
\end{array}$};
		\end{tikzpicture}

\end{minipage}
\begin{minipage}[t]{0.3\linewidth} 
\centering
\begin{tikzpicture}[xscale=0.8,yscale=0.8]
			
	\draw[blue,thick] (0,2)\nn to (0,2)\nn;
			
			\draw[thick,red] (0,0) to (2,1.5);	
			\draw[blue,thick] (0,-2)\nn to (0,-2)\nn;
				\draw[blue,thick] (0,0)\nn to (0,0)\nn;

\draw[thick,color=blue] (0,-2) .. controls (-2,-1) and (-2,1) .. (0,2);	
\draw[thick,color=blue] (0,-2) .. controls (2,-1) and (2,1) .. (0,2);	

\draw[thick,color=blue] (0,-2) .. controls (-0.5,-1.5) and (-0.5,-0.5) .. (0,0);	
\draw[thick,color=blue] (0,-2) .. controls (0.5,-1.5) and (0.5,-0.5) .. (0,0);	
		
\draw[blue](0,-2)node[below]{$P_1$}(0,2)node[above]{$P_2$}(0,0)node[left]{$P$};	
\draw[blue](0.2,-0.3)node[rotate=20]{$\bowtie$};	
\draw[red](0.4,0.3)node[rotate=-45]{$\bowtie$}(-0.4,0.3)node[rotate=45]{$\bowtie$};		

\draw[thick,color=red] (1.4,-1.8) .. controls (0.2,-1.3) and (-0.2,-1.3) .. (-1.4,-1.8);	
\draw[thick,color=red] (1,2) .. controls (0.5,1.5) and (-0.5,1.5) .. (-1,2);
	\draw[thick,red] (0,0) to (-2,1.5);	

\draw[blue](0.7,0.57)node{$l_{12}^\ast$}(0,1.6)node{$l_{03}^\ast$}(-0.7,0.57)node{$l_{11}^\ast$}(0,-1.4)node{$l_{04}^\ast$};
\node at (0,-4){(A5): $\begin{array}{c}
 l_{01}^\ast=l_{02}^\ast=0\\
l_{03}^\ast\neq 0, l_{09}^\ast= 0, l_{10}^\ast= 0
\end{array}$};

\end{tikzpicture}
\end{minipage}%
\caption{Case (A): $l_{01}^\ast+l_{02}^\ast+l_{03}^\ast\neq 0$}\label{f:Case A}
\end{figure}
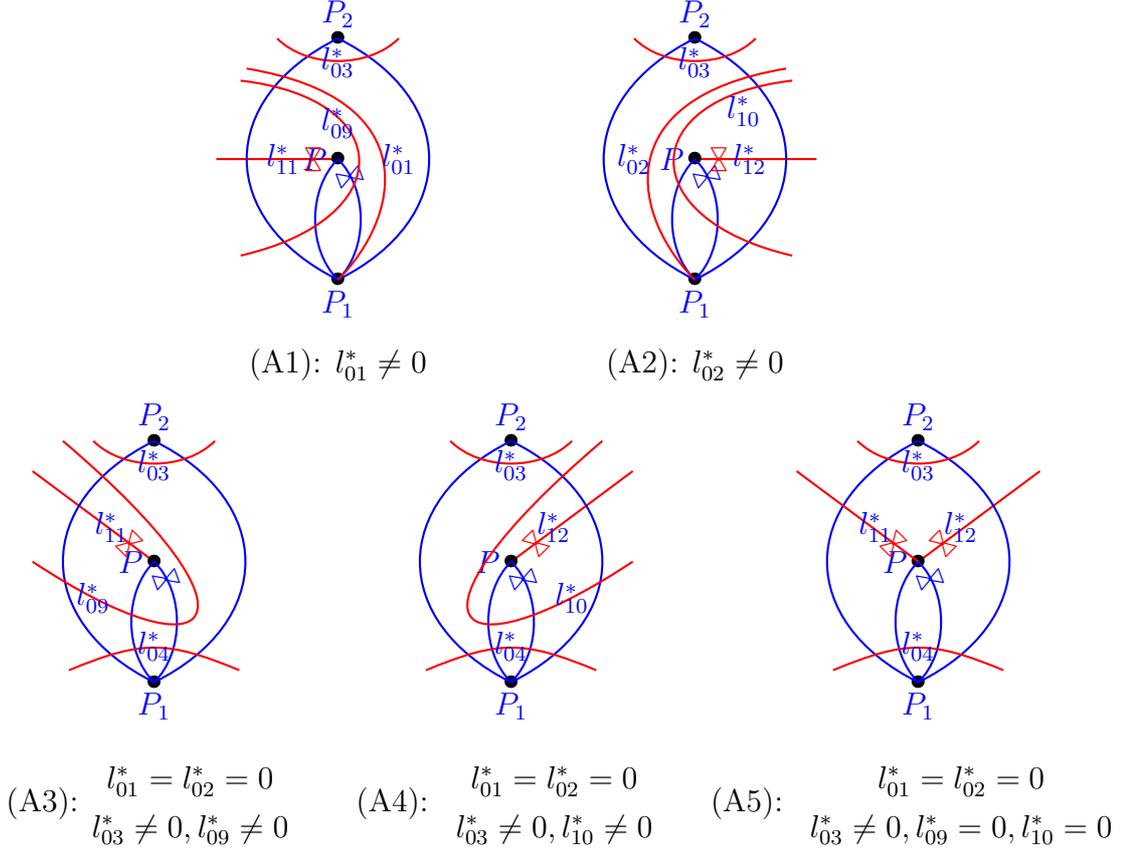
\begin{figure}[ht]
\begin{minipage}[t]{0.3\linewidth} 
\centering
\begin{tikzpicture}[xscale=0.8,yscale=0.8]
		\draw[blue,thick] (0,2)\nn to (0,2)\nn;
			
		\draw[thick,color=red] (0,2) .. controls (-0.8,1) and (-0.5,-0.3) .. (0,-0.3);						
\draw[thick,color=red] (0,2) .. controls (0.8,1) and (0.5,-0.3) .. (0,-0.3);	
\draw[blue](0.7,0.3)node{$l_{06}^\ast$};

		\draw[thick,red] (0,0)\nn to (0,2)\nn;
			\draw[blue,thick] (0,-2)\nn to (0,-2)\nn;
				\draw[blue,thick] (0,0)\nn to (0,0)\nn;

\draw[thick,color=blue] (0,-2) .. controls (-2,-1) and (-2,1) .. (0,2);	
\draw[thick,color=blue] (0,-2) .. controls (2,-1) and (2,1) .. (0,2);	

\draw[thick,color=blue] (0,-2) .. controls (-0.5,-1.5) and (-0.5,-0.5) .. (0,0);	
\draw[thick,color=blue] (0,-2) .. controls (0.5,-1.5) and (0.5,-0.5) .. (0,0);	
		
\draw[blue](0,-2)node[below]{$P_1$}(0,2)node[above]{$P_2$}(0,0)node[left]{$P$};	
\draw[blue](0.2,-0.3)node[rotate=20]{$\bowtie$};	
\draw[red](0,0.3)node[rotate=0]{$\bowtie$};		

\draw[thick,color=red] (1.4,-1.8) .. controls (0.2,-1.3) and (-0.2,-1.3) .. (-1.4,-1.8);	
\draw[thick,color=red] (2,-0.5) .. controls (-1,-2) and (-2,-1) .. (0,2);
\draw[blue](-1,-0.6)node{$l_{07}^\ast$}(0,1)node{$l_{05}^\ast$}(0,-1.4)node{$l_{04}^\ast$};
\node at (0,-3.4){(B1): $\begin{array}{c}l_{06}^\ast\neq 0,l_{07}^\ast\neq 0\end{array}$};
		\end{tikzpicture}
\end{minipage}%
\begin{minipage}[t]{0.3\linewidth} 
\centering
\begin{tikzpicture}[xscale=0.8,yscale=0.8]

			\draw[thick,color=red] (0,2) .. controls (-0.8,1) and (-0.5,-0.3) .. (0,-0.3);						
\draw[thick,color=red] (0,2) .. controls (0.8,1) and (0.5,-0.3) .. (0,-0.3);	
\draw[blue](-0.6,0.6)node{$l_{06}^\ast$};

		\draw[thick,red] (0,0)\nn to (0,2)\nn;
			\draw[blue,thick] (0,-2)\nn to (0,-2)\nn;
				\draw[blue,thick] (0,0)\nn to (0,0)\nn;
				
				\draw[blue,thick] (0,2)\nn to (0,2)\nn;
				
\draw[thick,color=blue] (0,-2) .. controls (-2,-1) and (-2,1) .. (0,2);	
\draw[thick,color=blue] (0,-2) .. controls (2,-1) and (2,1) .. (0,2);	

\draw[thick,color=blue] (0,-2) .. controls (-0.5,-1.5) and (-0.5,-0.5) .. (0,0);	
\draw[thick,color=blue] (0,-2) .. controls (0.5,-1.5) and (0.5,-0.5) .. (0,0);	
		
\draw[blue](0,-2)node[below]{$P_1$}(0,2)node[above]{$P_2$}(0,0)node[left]{$P$};	
\draw[blue](0.2,-0.3)node[rotate=20]{$\bowtie$};	
\draw[red](0,0.3)node[rotate=0]{$\bowtie$};

\draw[thick,color=red] (-1.4,-1.8) .. controls (-0.2,-1.3) and (0.2,-1.3) .. (1.4,-1.8);	

\draw[thick,color=red] (-2,-0.5) .. controls (1,-2) and (2,-1) .. (0,2);

\draw[blue](1,-0.6)node{$l_{08}^\ast$}(0,-1.4)node{$l_{04}^\ast$}(0,1)node{$l_{05}^\ast$};
\node at (0,-3.4){(B2): $\begin{array}{c}l_{06}^\ast\neq 0,l_{08}^\ast\neq 0\end{array}$};
\end{tikzpicture}

\end{minipage}
\begin{minipage}[t]{0.3\linewidth} 
\centering
\begin{tikzpicture}[xscale=0.8,yscale=0.8]
			
		\draw[thick,red] (0,0)\nn to (0,2)\nn;
		\draw[red](0,0.3)node[rotate=0]{$\bowtie$};
		
			\draw[blue,thick] (0,-2)\nn to (0,-2)\nn;

\draw[thick,color=red] (0,2) .. controls (-0.8,1) and (-0.5,-0.3) .. (0,-0.3);						
\draw[thick,color=red] (0,2) .. controls (0.8,1) and (0.5,-0.3) .. (0,-0.3);	
\draw[blue](0.7,0.3)node{$l_{06}^\ast$};

	\draw[thick,color=red] (1.4,-1.8) .. controls (0.2,-1.3) and (-0.2,-1.3) .. (-1.4,-1.8);				
\draw[thick,color=blue] (0,-2) .. controls (-2,-1) and (-2,1) .. (0,2);	
\draw[thick,color=blue] (0,-2) .. controls (2,-1) and (2,1) .. (0,2);	

\draw[thick,color=blue] (0,-2) .. controls (-0.5,-1.5) and (-0.5,-0.5) .. (0,0);	
\draw[thick,color=blue] (0,-2) .. controls (0.5,-1.5) and (0.5,-0.5) .. (0,0);	
\draw[blue](0,1)node{$l_{05}^\ast$}(0,-1.4)node{$l_{04}^\ast$};		
\draw[blue](0,-2)node[below]{$P_1$}(0,2)node[above]{$P_2$}(0,0)node[left]{$P$};	
	
\draw[blue](0.2,-0.3)node[rotate=20]{$\bowtie$};

\node at (0,-3.3){(B3): $\begin{array}{c}l_{06}^\ast\neq 0,l_{07}^\ast=l_{08}^\ast= 0\end{array}$};
		\end{tikzpicture}
\end{minipage}%

\begin{minipage}[t]{0.3\linewidth} 
\centering
\begin{tikzpicture}[xscale=0.8,yscale=0.8]
		\draw[blue,thick] (0,2)\nn to (0,2)\nn;
			
			\draw[thick,red] (0,0) to (2,1.5);	
		\draw[thick,red] (0,0)\nn to (0,2)\nn;
			\draw[blue,thick] (0,-2)\nn to (0,-2)\nn;
				\draw[blue,thick] (0,0)\nn to (0,0)\nn;

\draw[thick,color=blue] (0,-2) .. controls (-2,-1) and (-2,1) .. (0,2);	
\draw[thick,color=blue] (0,-2) .. controls (2,-1) and (2,1) .. (0,2);	

\draw[thick,color=blue] (0,-2) .. controls (-0.5,-1.5) and (-0.5,-0.5) .. (0,0);	
\draw[thick,color=blue] (0,-2) .. controls (0.5,-1.5) and (0.5,-0.5) .. (0,0);	
		
\draw[blue](0,-2)node[below]{$P_1$}(0,2)node[above]{$P_2$}(0,0)node[left]{$P$};	
\draw[blue](0.2,-0.3)node[rotate=20]{$\bowtie$};	
\draw[red](0.4,0.3)node[rotate=-45]{$\bowtie$};	
\draw[red](0,0.3)node[rotate=0]{$\bowtie$};		

\draw[thick,color=red] (1.4,-1.8) .. controls (0.2,-1.3) and (-0.2,-1.3) .. (-1.4,-1.8);	
\draw[thick,color=red] (2,-0.5) .. controls (-1,-2) and (-2,-1) .. (0,2);
\draw[blue](0.7,0.57)node{$l_{12}^\ast$}(-1,0.2)node{$l_{07}^\ast$}(0,1)node{$l_{05}^\ast$}(0,-1.4)node{$l_{04}^\ast$};
\node at (0,-4){(B4): $\begin{array}{c}
    l_{06}^\ast=0,\\
    l_{05}^\ast\neq 0\neq l_{07}^\ast.
\end{array}$};
		\end{tikzpicture}
\end{minipage}%
\begin{minipage}[t]{0.3\linewidth} 
\centering
\begin{tikzpicture}[xscale=0.8,yscale=0.8]

			\draw[thick,red] (0,0) to (-2,1.5);	
		\draw[thick,red] (0,0)\nn to (0,2)\nn;
			\draw[blue,thick] (0,-2)\nn to (0,-2)\nn;
				\draw[blue,thick] (0,0)\nn to (0,0)\nn;
				
				\draw[blue,thick] (0,2)\nn to (0,2)\nn;
				
\draw[thick,color=blue] (0,-2) .. controls (-2,-1) and (-2,1) .. (0,2);	
\draw[thick,color=blue] (0,-2) .. controls (2,-1) and (2,1) .. (0,2);	

\draw[thick,color=blue] (0,-2) .. controls (-0.5,-1.5) and (-0.5,-0.5) .. (0,0);	
\draw[thick,color=blue] (0,-2) .. controls (0.5,-1.5) and (0.5,-0.5) .. (0,0);	
		
\draw[blue](0,-2)node[below]{$P_1$}(0,2)node[above]{$P_2$}(0,0)node[left]{$P$};	
\draw[blue](0.2,-0.3)node[rotate=20]{$\bowtie$};	
\draw[red](-0.4,0.3)node[rotate=45]{$\bowtie$};	
\draw[red](0,0.3)node[rotate=0]{$\bowtie$};		

\draw[thick,color=red] (-1.4,-1.8) .. controls (-0.2,-1.3) and (0.2,-1.3) .. (1.4,-1.8);	
\draw[thick,color=red] (-2,-0.5) .. controls (1,-2) and (2,-1) .. (0,2);


\draw[blue](-0.7,0.57)node{$l_{11}^\ast$}(1,0.2)node{$l_{08}^\ast$}(0,-1.4)node{$l_{04}^\ast$}(0,1)node{$l_{05}^\ast$};
\node at (0,-4){(B5): $\begin{array}{c}
    l_{06}^\ast=0,\\
    l_{05}^\ast\neq 0\neq l_{08}^\ast.
\end{array}$};
\end{tikzpicture}

\end{minipage}
\begin{minipage}[t]{0.3\linewidth} 
\centering
\begin{tikzpicture}[xscale=0.8,yscale=0.8]
			
\draw[thick,red] (0,0) to (-2,1.5);	
				\draw[thick,red] (0,0) to (2,1.5);	
		\draw[thick,red] (0,0)\nn to (0,2)\nn;
			\draw[blue,thick] (0,-2)\nn to (0,-2)\nn;
				\draw[blue,thick] (0,0)\nn to (0,0)\nn;
				
	\draw[thick,color=red] (1.4,-1.8) .. controls (0.2,-1.3) and (-0.2,-1.3) .. (-1.4,-1.8);				
\draw[thick,color=blue] (0,-2) .. controls (-2,-1) and (-2,1) .. (0,2);	
\draw[thick,color=blue] (0,-2) .. controls (2,-1) and (2,1) .. (0,2);	

\draw[thick,color=blue] (0,-2) .. controls (-0.5,-1.5) and (-0.5,-0.5) .. (0,0);	
\draw[thick,color=blue] (0,-2) .. controls (0.5,-1.5) and (0.5,-0.5) .. (0,0);	
\draw[blue](-0.7,0.57)node{$l_{11}^\ast$}(0.7,0.57)node{$l_{12}^\ast$}(0,1)node{$l_{05}^\ast$}(0,-1.4)node{$l_{04}^\ast$};		
\draw[blue](0,-2)node[below]{$P_1$}(0,2)node[above]{$P_2$}(0,0)node[left]{$P$};	
\draw[blue](0.2,-0.3)node[rotate=20]{$\bowtie$};	
\draw[red](0,0.3)node{$\bowtie$}(0.4,0.3)node[rotate=-45]{$\bowtie$}(-0.4,0.3)node[rotate=45]{$\bowtie$};		

\node at (0,-4){(B6): $\begin{array}{c}
    l_{06}^\ast=0,l_{05}^\ast\neq 0,\\ l_{07}^\ast=l_{08}^\ast=0.
\end{array}$};
		\end{tikzpicture}
\end{minipage}%

\begin{minipage}[t]{0.4\linewidth} 
\centering
\begin{tikzpicture}[xscale=0.8,yscale=0.8]

			\draw[blue,thick] (0,2)\nn to (0,2)\nn;
			
			\draw[thick,red] (0,0) to (2,1.5);	
			\draw[blue,thick] (0,-2)\nn to (0,-2)\nn;
				\draw[blue,thick] (0,0)\nn to (0,0)\nn;

\draw[thick,color=blue] (0,-2) .. controls (-2,-1) and (-2,1) .. (0,2);	
\draw[thick,color=blue] (0,-2) .. controls (2,-1) and (2,1) .. (0,2);	

\draw[thick,color=blue] (0,-2) .. controls (-0.5,-1.5) and (-0.5,-0.5) .. (0,0);	
\draw[thick,color=blue] (0,-2) .. controls (0.5,-1.5) and (0.5,-0.5) .. (0,0);	
		
\draw[blue](0,-2)node[below]{$P_1$}(0,2)node[above]{$P_2$}(0,0)node[left]{$P$};	
\draw[blue](0.2,-0.3)node[rotate=20]{$\bowtie$};	
\draw[red](0.4,0.3)node[rotate=-45]{$\bowtie$};		

\draw[thick,color=red] (1.4,-1.8) .. controls (0.2,-1.3) and (-0.2,-1.3) .. (-1.4,-1.8);	
\draw[thick,color=red] (2,0) .. controls (-1,-2) and (-2,-1) .. (1.5,2);
\draw[thick,color=red] (2,-0.5) .. controls (-1,-2) and (-2,-1) .. (0,2);
\draw[blue](0.7,0.57)node{$l_{12}^\ast$}(-0.6,1)node{$l_{07}^\ast$}(1,-0.6)node{$l_{10}^\ast$}(0,-1.4)node{$l_{04}^\ast$};
\node at (0,-4.6){(B7): $\begin{array}{c}
l_{05}^\ast=l_{06}^\ast=0,\\
l_{07}^\ast\neq 0.
\end{array}
$ };
	\end{tikzpicture}
\end{minipage}%
\begin{minipage}[t]{0.4\linewidth} 
\centering
\begin{tikzpicture}[xscale=0.8,yscale=0.8]


			\draw[thick,red] (0,0) to (-2,1.5);	
			\draw[blue,thick] (0,-2)\nn to (0,-2)\nn;
				\draw[blue,thick] (0,0)\nn to (0,0)\nn;
				
				\draw[blue,thick] (0,2)\nn to (0,2)\nn;
				
\draw[thick,color=blue] (0,-2) .. controls (-2,-1) and (-2,1) .. (0,2);	
\draw[thick,color=blue] (0,-2) .. controls (2,-1) and (2,1) .. (0,2);	

\draw[thick,color=blue] (0,-2) .. controls (-0.5,-1.5) and (-0.5,-0.5) .. (0,0);	
\draw[thick,color=blue] (0,-2) .. controls (0.5,-1.5) and (0.5,-0.5) .. (0,0);	
		
\draw[blue](0,-2)node[below]{$P_1$}(0,2)node[above]{$P_2$}(0,0)node[left]{$P$};	
\draw[blue](0.2,-0.3)node[rotate=20]{$\bowtie$};	
\draw[red](-0.4,0.3)node[rotate=45]{$\bowtie$};		

\draw[thick,color=red] (-1.4,-1.8) .. controls (-0.2,-1.3) and (0.2,-1.3) .. (1.4,-1.8);	
\draw[thick,color=red] (-2,0) .. controls (1,-2) and (2,-1) .. (-1.5,2);
\draw[thick,color=red] (-2,-0.5) .. controls (1,-2) and (2,-1) .. (0,2);


\draw[blue](-0.7,0.57)node{$l_{11}^\ast$}(0.6,1)node{$l_{08}^\ast$}(-0.9,-0.5)node{$l_{09}^\ast$}(0,-1.4)node{$l_{04}^\ast$};
\node at (0,-4.6){(B8): $\begin{array}{c}
l_{05}^\ast=l_{06}^\ast=0,\\
l_{08}^\ast\neq 0.
\end{array}
$ };

	\end{tikzpicture}
\end{minipage}%

\caption{Case (B): $\begin{cases}
l_{01}^\ast+l_{02}^\ast+l_{03}^\ast=0\\
 l_{05}^\ast+l_{06}^\ast+l_{07}^\ast+l_{08}^\ast\neq 0
\end{cases}$
}\label{f:Case B}
\end{figure}
\begin{figure}[ht]
\begin{minipage}[t]{0.3\linewidth} 
\centering
\begin{tikzpicture}[xscale=0.8,yscale=0.8]


			\draw[thick,red] (0,0) to (-2,1.5);	
			\draw[blue,thick] (0,-2)\nn to (0,-2)\nn;
				\draw[blue,thick] (0,0)\nn to (0,0)\nn;
				
				\draw[blue,thick] (0,2)\nn to (0,2)\nn;
				
\draw[thick,color=blue] (0,-2) .. controls (-2,-1) and (-2,1) .. (0,2);	
\draw[thick,color=blue] (0,-2) .. controls (2,-1) and (2,1) .. (0,2);	

\draw[thick,color=blue] (0,-2) .. controls (-0.5,-1.5) and (-0.5,-0.5) .. (0,0);	
\draw[thick,color=blue] (0,-2) .. controls (0.5,-1.5) and (0.5,-0.5) .. (0,0);	
		
\draw[blue](0,-2)node[below]{$P_1$}(0,2)node[above]{$P_2$}(0,0)node[left]{$P$};	
\draw[blue](0.2,-0.3)node[rotate=20]{$\bowtie$};	
\draw[red](-0.4,0.3)node[rotate=45]{$\bowtie$};		

\draw[thick,color=red] (-1.4,-1.8) .. controls (-0.2,-1.3) and (0.2,-1.3) .. (1.4,-1.8);	
\draw[thick,color=red] (-2,0) .. controls (1,-2) and (2,-1) .. (-1.5,2);


\draw[blue](-0.7,0.57)node{$l_{11}^\ast$}(-0.9,-0.5)node{$l_{09}^\ast$}(0,-1.4)node{$l_{04}^\ast$};
\node at (0,-3.4){(C1): $l_{09}^\ast\neq 0
$ };

	\end{tikzpicture}
\end{minipage}%
\begin{minipage}[t]{0.3\linewidth} 
\centering
\begin{tikzpicture}[xscale=0.8,yscale=0.8]

			\draw[blue,thick] (0,2)\nn to (0,2)\nn;
			
			\draw[thick,red] (0,0) to (2,1.5);	
			\draw[blue,thick] (0,-2)\nn to (0,-2)\nn;
				\draw[blue,thick] (0,0)\nn to (0,0)\nn;

\draw[thick,color=blue] (0,-2) .. controls (-2,-1) and (-2,1) .. (0,2);	
\draw[thick,color=blue] (0,-2) .. controls (2,-1) and (2,1) .. (0,2);	

\draw[thick,color=blue] (0,-2) .. controls (-0.5,-1.5) and (-0.5,-0.5) .. (0,0);	
\draw[thick,color=blue] (0,-2) .. controls (0.5,-1.5) and (0.5,-0.5) .. (0,0);	
		
\draw[blue](0,-2)node[below]{$P_1$}(0,2)node[above]{$P_2$}(0,0)node[left]{$P$};	
\draw[blue](0.2,-0.3)node[rotate=20]{$\bowtie$};	
\draw[red](0.4,0.3)node[rotate=-45]{$\bowtie$};		

\draw[thick,color=red] (1.4,-1.8) .. controls (0.2,-1.3) and (-0.2,-1.3) .. (-1.4,-1.8);	
\draw[thick,color=red] (2,0) .. controls (-1,-2) and (-2,-1) .. (1.5,2);
\draw[blue](0.7,0.57)node{$l_{12}^\ast$}(1,-0.6)node{$l_{10}^\ast$}(0,-1.4)node{$l_{04}^\ast$};
\node at (0,-3.4){(C2): $l_{10}^\ast\neq 0$ };
	\end{tikzpicture}
\end{minipage}%
\begin{minipage}[t]{0.3\linewidth} 
\centering
\begin{tikzpicture}[xscale=0.8,yscale=0.8]
			
	\draw[blue,thick] (0,2)\nn to (0,2)\nn;
			
			\draw[thick,red] (0,0) to (2,1.5);	
			\draw[blue,thick] (0,-2)\nn to (0,-2)\nn;
				\draw[blue,thick] (0,0)\nn to (0,0)\nn;

\draw[thick,color=blue] (0,-2) .. controls (-2,-1) and (-2,1) .. (0,2);	
\draw[thick,color=blue] (0,-2) .. controls (2,-1) and (2,1) .. (0,2);	

\draw[thick,color=blue] (0,-2) .. controls (-0.5,-1.5) and (-0.5,-0.5) .. (0,0);	
\draw[thick,color=blue] (0,-2) .. controls (0.5,-1.5) and (0.5,-0.5) .. (0,0);	
		
\draw[blue](0,-2)node[below]{$P_1$}(0,2)node[above]{$P_2$}(0,0)node[left]{$P$};	
\draw[blue](0.2,-0.3)node[rotate=20]{$\bowtie$};	
\draw[red](0.4,0.3)node[rotate=-45]{$\bowtie$}(-0.4,0.3)node[rotate=45]{$\bowtie$};		

\draw[thick,color=red] (1.4,-1.8) .. controls (0.2,-1.3) and (-0.2,-1.3) .. (-1.4,-1.8);	
	\draw[thick,red] (0,0) to (-2,1.5);	

\draw[blue](0.7,0.57)node{$l_{12}^\ast$}(-0.7,0.57)node{$l_{11}^\ast$}(0,-1.4)node{$l_{04}^\ast$};
\node at (0,-3.4){(C3): $l_{09}^\ast=l_{10}^\ast=0$};

\end{tikzpicture}
\end{minipage}%
\caption{Case (C): $-l_{04}^\ast+\sum_{i=1}^8l_{0i}^\ast=0$}\label{f:Case C}
\end{figure}
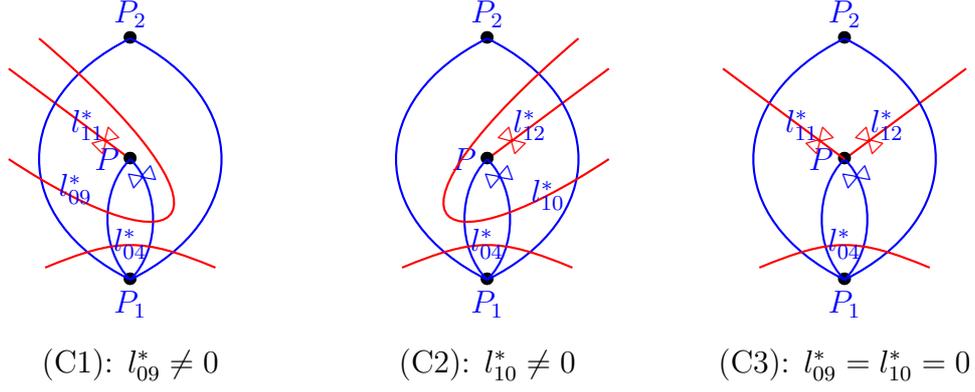
\begin{lem}\label{lem:correspodence-type}
The following statements hold:
\begin{enumerate}
    \item $\Arcs_\Delta(\cm)$ satisfies  $({\rm A}i)$ if and only if $\Arcs_\Delta(\cn)$ satisfies $({\rm A}i)$ for $1\leq i\leq 5$;
    \item $\Arcs_\Delta(\cm)$ satisfies  $({\rm B}i)$ if and only if $\Arcs_\Delta(\cn)$ satisfies $({\rm B}i)$ for $1\leq i\leq 8$;
    \item $\Arcs_\Delta(\cm)$ satisfies $({\rm C}i)$ if and only if $\Arcs_\Delta(\cn)$ satisfies $({\rm C}i)$ for $1\leq i\leq 3$.
\end{enumerate}
\end{lem}
\begin{proof}
 For $i,j,k,l\in \{-1,0,1\}$, we define
\[
D_{i,j,k,l}(\ast):=i\times \Int(\ba_1|\ast)+j\times \Int(\ba_2|\ast)+k\times\Int(\ba_3|\ast)+l\times\Int(\ba_4|\ast),
\]
 where $\ast\in \{\cm,\cn\}$.
A direct computation shows that
\begin{align}\label{eq:-1-234}
    D_{-1,-1,1,1}(\ast)=\begin{cases}
        -2l_{03}^\ast-l_{01}^\ast-l_{02}^\ast<0&\text{if $\ast$ satisfies one of $({\rm A})$};\\
       l_{05}^\ast+2l_{06}^\ast+l_{07}^\ast+l_{08}^\ast >0&\text{if $\ast$ satisfies one of $({\rm B})$};\\
        0&\text{if $\ast$ satisfies one of $({\rm C})$}.
    \end{cases}
\end{align}

Similarly, we can compute
\begin{align}
    D_{-1,1,1,1}(\ast)&=\begin{cases}\label{eq:1-234}
        -l_{01}^\ast<0 &\text{case $({\rm A}1)$};\\
        l_{02}^\ast+4l_{10}^\ast+2l_{12}^\ast>0 &\text{case $({\rm A}2)$};\\
        2l_{04}^\ast\geq 0&\text{case $({\rm A}3)$};\\
       2l_{04}^\ast+4l_{10}^\ast+2l_{12}^\ast>0&\text{case $({\rm A}4)$};\\
       2l_{12}^\ast+2l_{04}^\ast\geq 0&\text{case $({\rm A}5)$}.
    \end{cases}\\
D_{1,-1,1,-1}(\ast)&=\begin{cases}\label{eq:-12}
    -l_{02}^\ast-2l_{10}^\ast<0 &\text{case $({\rm A}2)$};\\
      2l_{09}^\ast+2l_{11}^\ast>0 &\text{case $({\rm A}3)$};\\
    -2l_{10}^\ast<0&\text{case $({\rm A}4)$};\\
    2l_{11}^\ast\geq 0&\text{case $({\rm A}5)$}.
\end{cases}\\
    D_{1,-1,-1,1}(\ast)&=\begin{cases}\label{eq:1-2-34}
    2l_{09}^\ast>0 &\text{case $({\rm A}3)$};\\
    -2l_{12}^\ast\leq 0& \text{case $({\rm A}5)$}.
\end{cases}\\
    D_{-1,1,-1,-1}(\ast)&=\begin{cases}\label{eq:-12-3-4}
 l_{02}^\ast>0 &\text{case $({\rm A}2)$};\\
-2l_{04}^\ast\leq 0&\text{case $({\rm A}4)$}.
\end{cases}
\end{align}

By $(\ref{eq:-1-234})$-$(\ref{eq:-12-3-4})$, we conclude that $\Arcs_\Delta(\cm)$ satisfies $({\rm A}i)$ if and only if $\Arcs_\Delta(\cn)$ satisfies $({\rm A}i)$ for $1\leq i\leq 5$. 

Now we turn to the case $({\rm B})$. A direct computation shows that
\begin{align}
    D_{-1,0,1,0}(\ast)&=\begin{cases}l_{05}^\ast+l_{06}^\ast+l_{07}^\ast+l_{10}^\ast+l_{12}^\ast>0& \text{cases $({\rm B}1)$-$({\rm B}7)$;}\\
    -l_{09}^\ast\leq 0 &\text{case $({\rm B}8)$}.
    \end{cases}\label{eq:-1010}\\
    D_{0,-1,1,0}(\ast)&=\begin{cases}
        l_{05}^\ast+l_{06}^\ast+l_{08}^\ast+l_{11}^{\ast}>0& \text{cases $({\rm B}1)$-$({\rm B}6)$};\\
        -l_{10}^\ast\leq 0&\text{case $({\rm B}7)$}.
    \end{cases}\label{eq:0-110}\\
    D_{0,-1,0,1}(\ast)&=\begin{cases}
        l_{06}^\ast+l_{08}^\ast>0& \text{cases $({\rm B}1)$-$({\rm B}3)$ and $({\rm B}5)$};\\
        -l_{12}^\ast\leq 0 &\text{cases $({\rm B}4)$ and $({\rm B}6)$}.
    \end{cases}\label{eq:0-101}\\
    D_{1,0,0,-1}(\ast)&=\begin{cases}
        -l_{07}^\ast<0 &\text{case $({\rm B}4)$};\\
        l_{11}^\ast \geq 0 &\text{case $({\rm B}6)$}.
    \end{cases}\label{eq:100-1}\\
    D_{-1,0,0,1}(\ast)&=\begin{cases}
        l_{06}^\ast+l_{07}^\ast>0 &\text{cases $({\rm B}1)$-$({\rm B}3)$};\\
        -l_{11}^\ast\leq 0 &\text{case $({\rm B}5)$}.
    \end{cases}\label{eq:-1001}\\
    D_{1,-1,0,0}(\ast)&=\begin{cases}
        -l_{07}^\ast<0&\text{case $({\rm B}1)$};\\
        l_{08}^\ast>0 &\text{case $({\rm B}2)$};\\
        0&\text{case $({\rm B}3)$}.
    \end{cases}\label{eq:1-100}
\end{align}
By (\ref{eq:-1-234}) and (\ref{eq:-1010})-(\ref{eq:1-100}), we conclude that $\Arcs_\Delta(\cm)$ satisfies $({\rm B}i)$ if and only if $\Arcs_\Delta(\cn)$ satisfies $({\rm B}i)$ for $1\leq i\leq 8$. 

Finally, consider the case $({\rm C})$. We have
\begin{align}
    D_{1,0,-1,0}(\ast)&=\begin{cases}
    l_{09}^\ast>0&\text{case $({\rm C}1)$};\\
    -l_{10}^\ast-l_{12}^\ast<0&\text{case $({\rm C}2)$};\\
    -l_{12}^\ast\leq 0&\text{case $({\rm C}3)$}.
    \end{cases}\label{eq:10-10}\\
    D_{0,1,-1,0}(\ast)&=\begin{cases}
        l_{10}^\ast>0&\text{case $({\rm C}2)$};\\
        -l_{11}^\ast\leq 0&\text{case $({\rm C}3)$}.
    \end{cases}\label{eq:01-10}
\end{align}
According to (\ref{eq:-1-234}), (\ref{eq:10-10}) and (\ref{eq:01-10}), we obtain that $\Arcs_\Delta(\cm)$ satisfies $({\rm C}i)$ if and only if $\Arcs_\Delta(\cn)$ satisfies $({\rm C}i)$ for $1\leq i\leq 3$.  This completes the proof.  
\end{proof}

\begin{prop}\label{prop:1-puncture}
 Let $\Delta$ be a $1$-puncture piece, then $\Arcseg_{\Delta}(\cm)=\Arcseg_{\Delta}(\cn)$.
\end{prop}
\begin{proof} 
Keep the notation as above. It suffices to show that $l_{ij}^\cm=l_{ij}^\cn$ for each $ij$.
This can be verified case by case, according to Table \ref{table:16-cases}. Let us  give the verification for the case $({\rm A}1)$ and leave the other tedious computations as an exercise.
Let $d_i:=\Int(\ba_i|\cm)=\Int(\ba_i|\cn)$ for $1\leq i\leq 4$.
Assume that $\Arcs_\Delta(\cm)$ satisfies $({\rm A}1)$. By Lemma \ref{lem:correspodence-type}, $\Arcs_\Delta(\cn)$ also satisfies $({\rm A}1)$. According to the definition of intersection number, we have
\begin{equation}\label{eq:A1}
    \begin{cases}
        d_1=l_{01}^\ast+l_{03}^\ast+2l_{09}^\ast+l_{11}^\ast;\\
        d_2=l_{03}^\ast;\\
        d_3=l_{09}^\ast+l_{11}^\ast;\\
        d_4=l_{09}^\ast,
    \end{cases}
\end{equation}
where $\ast\in \{\cm,\cn\}$.
By considering $l_{01}^\ast$,$l_{03}^\ast,l_{09}^\ast,l_{11}^\ast$ as variables, it is easy to see that the equation (\ref{eq:A1}) has a unique solution. Consequently, $l_{01}^\cm=l_{01}^\cn$, $l_{03}^\cm=l_{03}^\cn$, $l_{09}^\cm=l_{09}^\cn$ and $l_{11}^\cm=l_{11}^\cn$. Note that all the other $l_{ij}^\ast$ that do not appear in the case $({\rm A}1)$ of Figure \ref{f:Case A} are zeros, since $\cm$ and $\cn$ are permissible. This completes the verification for the case $({\rm A}1)$.
\end{proof}

\section{Intersection vector determines arcs}\label{s:inter-vector-determine-arc}
\subsection{Disc with at most one puncture}
Recall that $(\fS,\fM)$ is a disc with at most one puncture $P$ and $\fT$ is a tagged triangulation. In particular, $\fT$ contains no loops. Hence for any tagged arc $\ba\in \fT$ and $\beta \in \mathbb{A}_{\bowtie}(\fS)$, we have
\[\Int(\alpha|\beta)=\Int^A(\alpha|\beta)+\Int^C(\alpha|\beta)+\Int^D(\alpha|\beta).\]

\begin{lem}\label{l: T has no element in cm and cn} Let $\cm,\cn$ be two finite multisets consisting of pairwise compatible admissible tagged arcs such that $\Intv_{\TT}(\cm)=\Intv_{\TT}(\cn)$. 
Let $\ba \in \TT$ be an arc.  Assume that $\ba^{(m)} \subseteq \cm, \ba^{(m+1)}\not\subseteq\cm$ and $\ba^{(n)} \subseteq \cn, \ba^{(n+1)}\not\subseteq\cn$ for some non-negative integers $m$ and $n$, then $m = n$.

\end{lem}
\begin{proof} 
Without loss of generality, we may assume that $m>0$. We are going to show that $m=n$.
Given that all tagged arcs  in $\cm$ are compatible,
   we have $\Int(\ba|\gamma)=0$ for any $\gamma\in \cm$ with $\gamma\neq \ba$. On the other hand, $\Int(\ba|\ba)=-1$.
It follows that $\Int(\ba|\cm)=-m$. 
Since $\Intv_{\TT}(\cm)=\Intv_{\TT}(\cn)$, we conclude that $\Int(\ba|\cn)=-m$.  We claim that $n>0$. Otherwise $n=0$. In other words, $\ba\not\in \cn$. According to the definition of intersection number, $\Int(\ba|\cn)\geq 0$, which contradicts the fact that $\dInt(\ba|\cn)=-m<0$. Hence $n>0$.
A similar discussion yields $\dInt(\ba|\cn)=-n$, and hence $m=n$.
\end{proof}

\begin{thm}\label{t: main theorem}
    Let $\cm,\cn$ be two finite multisets consisting of pairwise compatible admissible tagged arcs such that $\dInttv_{\TT}(\cm)=\dInttv_{\TT}(\cn)$, then $\cm=\cn$.
\end{thm}
\begin{proof}
 
By Lemma \ref{l: T has no element in cm and cn}, we may assume that $\cm\cap \TT=\emptyset=\cn\cap\TT$.
Recall that $\cm^\ast$ (resp. $\cn^\ast$) is the multiset obtained from $\cm$ (resp. $\cn$) by replacing each conjugate pair $(\gamma_Q^-,\gamma_Q^{\bowtie})$ in $\cm$ (resp. $\cn$) by the corresponding loop $l_{Q}$ which closely wraps around $\gamma_Q^-$. Clearly, $\cm^\ast$ and $\cn^\ast$ are permissible with respect to $\fT$ (cf. Section \ref{l: T has no element in cm and cn}).
We claim that $\cm^*=\cn^*$.

According to Lemma~\ref{l: m* has same intersection vector with m}, we have
 \[\dInttv_{\TT}(\cm^{*})=\dInttv_{\TT}(\cm)= \dInttv_{\TT}(\cn)=\dInttv_{\TT}(\cn^{*}).\]
 It follows that $\Arcs_{\Delta}(\cm^*)=\Arcs_{\Delta}(\cn^*)$
for each tile $\Delta$ of $(\fS,\fM,\fT)$ by Theorem \ref{thm:intersection-arc-segments}.
By induction on $|\cm^\ast|$, it suffices to prove that for any $\alpha\in \cm^*$, there is an arc $\beta\in \cn^*$ such that $\alpha=\beta$.  

Let us assume first that $\cm^\ast$ contains arcs which have two different boundary endpoints. Since $\cm^\ast$ is finite and $\fS$ is a disc, there exists an arc $\alpha\in \cm^\ast$ with a fixed orientation such that no endpoints of arcs in  $\cm^\ast$ lie on the right-hand side of $\alpha$. It follows that there are no endpoints of arcs in $\cn^\ast$ lying on the right hand side of $\alpha$ by the fact that $\Arcs_\Delta(\cm^\ast)=\Arcs_\Delta(\cn^\ast)$ for any tile $\Delta$.  Consequently, $\alpha$ is compatible with any arc in $\cn^\ast$.

Denote by $\alpha(0)=P_0$ and $\alpha(1)=Q$. Let $P_1,\dots, P_m$ be the intersections of $\alpha$ with $\fT$ in order, which belong to $\ba_1,\dots, \ba_m\in \fT$ respectively. We also denote by $\alpha_{(1)},\dots, \alpha_{(m+1)}$ the irreducible arc segments of $\alpha$ in order, which belong to the tile $\Delta_1,\dots, \Delta_{m+1}$, respectively (cf. Figure \ref{f:exist of beta}). Since $\Arcs_{\Delta_1}(\cm^\ast)=\Arcs_{\Delta_1}(\cn^\ast)$, there are arcs in $\cn^\ast$ that admit $\alpha_{(1)}$ as an irreducible arc segment. Since $\cn^\ast$ is finite, we may choose an arc $\beta\in \cn^\ast$ such that $\beta$ has the maximal number of common consecutive irreducible arc segments with $\alpha$ in order. Let $\beta_{(1)}, \dots, \beta_{m'+1}$ be all the irreducible arc segments of $\beta$ in order such that $\beta_{(1)}=\alpha_{(1)}$.  If $m'=m$ and $\beta_{(i)}=\alpha_{(i)}$ for all $1\leq i\leq m+1$, then $\beta=\alpha$. Otherwise, $\beta_{(1)}=\alpha_{(1)}, \dots,\beta_{(l)}=\alpha_{(l)}$ and $\beta_{l+1}\neq \alpha_{(l+1)}$ for some positive integer $l$. Note that $\Arcs_{\Delta_{l+1}}(\cm^\ast)=\Arcs_{\Delta_{l+1}}(\cn^\ast)$, there is an arc $\beta'\in \cn^\ast$ which admits $\alpha_{(l+1)}$ as an irreducible arc segment. Since $\alpha$ is compatible with any arc in $\cn^\ast$ and any two arcs in $\cn^\ast$ are compatible, we conclude that $\beta'$ has at least $l+1$ common consecutive irreducible arc segments with $\alpha$ in order (cf. Figure \ref{f:exist of beta}), which contradicts the maximality of $\beta$. Hence $\beta=\alpha$.

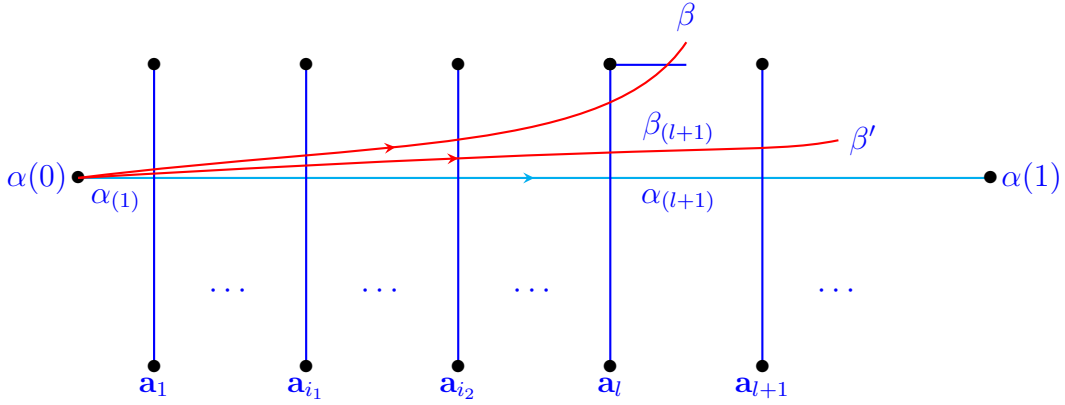
\begin{figure}[ht]
		\begin{tikzpicture}
			\draw[cyan,thick,->-=.5,>=stealth] (-6,0.5)\nn to (6,0.5)\nn;
			\draw[blue] (-6,0.5)node[left]{$\alpha(0)$} (6,0.5)node[right]{$\alpha(1)$};
				\draw[blue,thick] (-5,-2)\nn to (-5,2)\nn;
			\draw[blue,thick] (-3,-2)\nn to (-3,2)\nn;
			

\draw[blue,thick] (-1,-2)\nn to (-1,2)\nn;
			
\draw[blue,thick] (3,-2)\nn to (3,2)\nn;

			\draw[blue] (0,-1)node{$\cdots$}  (-4,-1)node{$\cdots$}  (-2,-1)node{$\cdots$}    (4,-1)node{$\cdots$} ;
			\draw[blue,thick] (1,-2)\nn to (1,2)\nn;
				\draw[blue,thick] (1,2)\nn to (2,2);
			\draw[blue] (-5,-2)node[below]{$\ba_{1}$} (-3,-2)node[below]{$\ba_{i_1}$}  (-1,-2)node[below]{$\ba_{i_2}$}   (1,-2)node[below]{$\ba_l$} (3,-2)node[below]{$\ba_{l+1}$};
			\draw[blue] (-5.5,0.5)node[below]{$\alpha_{(1)}$}  (1.9,0.8)node[above]{$\beta_{(l+1)}$} (1.9,0.5)node[below]{$\alpha_{(l+1)}$};
		
			\draw[thick,red,->-=.5,>=stealth] (-6,0.5) .. controls (-2,1) and (1,0.8) .. (2,2.3);
			\draw[thick,red,->-=.5,>=stealth] (-6,0.5) .. controls (2.5,1) and (3,0.8) .. (4,1);
			
			\draw[blue](2,2.3)node[above]{$\beta$}(4,1)node[right]{$\beta'$};


		\end{tikzpicture}
		\caption{Existence of $\beta$}\label{f:exist of beta}
	\end{figure}

 It remains to assume that each arc in $\cm^\ast$ and $\cn^\ast$ is either a once-punctured monogon with a boundary marked point as the base point, or an arc with one endpoint on the boundary and the other as the puncture $P$. If $\cm^\ast$ contains an arc $\alpha$ with $P$ and a boundary marked point $Q$ as endpoints, then by $\Arcs_\Delta(\cm^\ast)=\Arcs_\Delta(\cn^\ast)$, there are arcs in $\cn^\ast$ with $Q$  as an endpoint. If there is an arc $\beta\in \cn^\ast$ with $P,Q$ as endpoints, it turns out that $\beta=\alpha$. Otherwise, 
 all arcs in $\cn^\ast$ are once-punctured loops with $Q$ as the base point, since all arcs in $\cn^\ast$ are pairwise compatible. It follows that there are no arcs in $\cn^\ast$ with $P$ as an endpoint, a contradiction.  

 Now assume that $\cm^\ast$ and $\cn^\ast$ consist of once-punctured loops with boundary marked points as base points respectively. Since $\cm^\ast$ and $\cn^\ast$ are permissible and $\Arcs_\Delta(\cm^\ast)=\Arcs_\Delta(\cn^\ast)$, we conclude that $\cm^\ast=\alpha^{(m)}$ and $\cn^\ast=\alpha^{(m')}$ for some non-negative integers $m$ and $m'$, where $\alpha$ is a once-punctured loop with a boundary marked point as the base point. It is easy to see that $m=m'$ by considering the intersection number with arcs in $\fT$. Namely, for any $\alpha\in \cm^\ast$, there  is an arc $\beta\in \cn^\ast$ such that $\alpha=\beta$. This completes the proof.
 \end{proof}

\subsection{Disc with one unmarked boundary component in interior}
In this subsection, we generalize certain definitions and results to the setting of disc with one unmarked boundary component in its interior. 
Let $\fS$ be a disc with one unmarked boundary component in its interior, $\fM$ a non-empty finite set of marked points on the boundary $\partial \fS$. We assume that $|\fM|\geq 3$. The notions of {\em loop}, {\em intersection number} $\Int^\circ$, {\em compatible} and {\em arc} are  defined in the same way as  Section \ref{ss:disc-one-puncture}.

A triangulation $\fT$ of $(\fS,\fM)$ is a maximal set of pairwise compatible arcs. In this case, $\fS$ is divided by $\fT$ into a collection of the following tiles:
\begin{itemize}
    \item[(1)] monogons, that is, loops, with exactly one unmarked boundary component in their interior;
    \item[(2)] triangles.
\end{itemize}
Let $\fT$ be a triangulation of $(\fS,\fM)$. An arc is {\em permissible} with respect to $\fT$ if $\gamma\not\in \fT$.  Each permissible arc is divided by $\fT$ into irreducible arc segments.
For a tile $\Delta$ and a finite multiset $\cm$ of permissible arcs, we also denote by $\Arcs_{\Delta}(\cm)$ the finite multiset of all the irreducible arc segments of arcs in $\cm$ lying in $\Delta.$
The {\em intersection vector} of $\cm$ with respect to $\fT$ is defined as
\[
\Intv_\fT^\circ(\cm):=(\sum_{\gamma\in \cm}\Int^\circ(\ba|\gamma))_{\ba\in \fT}.
\]
\begin{prop}
    Let $\cm$ and $\cn$ be two finite multisets consisting of pairwise compatible permissible arcs respectively. Assume that $\Intv_\fT^\circ(\cm)=\Intv_\fT^\circ(\cn)$. Then $\Arcs_{\Delta}(\cm)=\Arcs_\Delta(\cn)$ for each tile $\Delta$.
\end{prop}
\begin{proof}
    If $\Delta$ is a triangle, then the result follows from the proof of Lemma \ref{l:tri-gons}. Suppose that  $\Delta$ is a monogon with loop $\ba$ and denote by $\angle\eta$ the unique angle. By the convention for arcs, each arc has no self-intersection, we conclude that $\Arcs_{\Delta}(\cm)$ consists of irreducible arc segments which cross the angle $\angle\eta$. It follows that 
    \[2|\Arcs_{\Delta}(\cm)|=\Int^\circ(\ba|\cm)=\Int^\circ(\ba|\cn)=2|\Arcs_{\Delta}(\cn)|.\]
    Hence, $\Arcs_{\Delta}(\cm)=\Arcs_{\Delta}(\cn)$. 
\end{proof}
A simplified  proof of Theorem \ref{t: main theorem} yields the following.
\begin{thm}\label{thm:type-c-intersection-vector}
     Let $\cm$ and $\cn$ be two finite multisets consisting of pairwise compatible permissible arcs respectively. Assume that $\Intv_\fT^\circ(\cm)=\Intv_\fT^\circ(\cn)$. Then $\cm=\cn$.
\end{thm}

\section{Proof of the main result}

\subsection{Fomin--Shapiro--Thurston correspondence}\label{ss:FST-corr}
Let $(\fS,\fM,\fP)$ be a marked surface with non-empty boundary, $\fT$ a tagged triangulation. Fomin, Shapiro and Thurston \cite{FST08} constructed a cluster algebra $\mathcal{A}_\fT$ associated with $(\fS,\fM,\fT)$. We refer to \cite{FST08} for the precise construction and  \cite{FZCA4} for unexplained terminology in cluster algebras. 
Moreover,
\begin{itemize}
    \item If $\fS$ is a disc, $|\fM|=n+3$ and $|\fP|=0$, then $\mathcal{A}_\fT$ is of type $\mathbb{A}_n$. 
    \item If $\fS$ is a disc, $|\fM|=n\geq 4$ and $|\fP|=1$, then $\mathcal{A}_\fT$ is of type $\mathbb{D}_n$.
\end{itemize}
The following is the so-called Fomin--Shapiro--Thurston correspondence.
\begin{thm}\cite[Theorem 7.11 \& Thoerem 8.6]{FST08}\label{thm:FST-correspondence}
    There is a bijection
    \[
    x_\fT[-]:\{\text{admissible tagged arcs in $\fS$}\}\to \{\text{cluster variables of $\mathcal{A}_\fT$}\},    \]
   which induces bijections
   \[
   x_\fT[-]:\{\text{tagged triangulations of $\fS$}\}\to \{\text{cluster of $\mathcal{A}_\fT$}\},
   \]
    \[
    x_\fT[-]:\begin{Bmatrix}\begin{array}{c}\text{finite multisets of pairwise compatible}\\ \text{admissible tagged arcs in $\fS$}\end{array}\end{Bmatrix}\to \begin{Bmatrix}\text{cluster monomials of $\mathcal{A}_\fT$}\end{Bmatrix}.
    \]
    Furthermore, let $\cm$ be a finite multiset consisting of pairwise compatible admissible tagged arcs in $\fS$, then the denominator vector of $x_\fT[\cm]$ with respect to the initial cluster $x_\fT[\fT]$ coincides with $\Intv_\fT(\cm)$.
\end{thm}




Now assume that $\fS$ is a disc with $n+1$ boundary marked points and an unmarked boundary component in its interior. Let $\fT$ be a triangulation of $(\fS,\fM)$. Denote by  $\fT=\{\ba_1,\ba_2,\dots,\ba_n\}$, where $\ba_1$ is the unique loop in $\fT$. Similar to \cite{FST08}, we can associated a cluster algebra $\mathcal{A}_\fT$ to $(\fS,\fM,\fT)$. The resulting cluster algebra is of type $\mathbb{C}_n$ (cf. \cite{FG22}). We have a similar Fomin--Shapiro--Thurston correspondence in this case, which is a consequence of \cite[Proposition 5.11]{FG22} and \cite[Theorem B]{HZZ23}.
\begin{prop}\label{prop:type-c-FST-correspondence}
    Let $\fS$ be a disc with $n+1$ boundary marked points and an unmarked boundary component in its interior. Denote by $\fT=\{\ba_1,\dots,\ba_n\}$ a triangulation, where $\ba_1$ is a loop. There are bijections
    \[
    x_\fT[-]:\{\text{arcs in $\fS$}\}\to \{\text{cluster variables of $\mathcal{A}_\fT$}\},
    \]
    \[
   x_\fT[-]:\{\text{triangulations of $\fS$}\}\to \{\text{cluster of $\mathcal{A}_\fT$}\},
   \]
    \[
    x_\fT[-]:\begin{Bmatrix}\begin{array}{c}\text{finite multisets of pairwise }\\ \text{ compatible arcs in $\fS$}\end{array}\end{Bmatrix}\to \begin{Bmatrix}\text{cluster monomials of $\mathcal{A}_\fT$}\end{Bmatrix}.
    \]
    Furthermore,  let $\cm$ be a finite multiset consisting of pairwise compatible permissible arcs in $\fS$ with intersection vector $\Intv_\fT^\circ(\cm)=(m_1,\dots, m_n)$, then the denominator vector of $x_\fT[\cm]$ with respect to the initial cluster $x_\fT[\fT]$ is $(\frac{m_1}{2},m_2,\dots, m_n)$.
\end{prop}

\subsection{Proof of Theorem \ref{thm:denominator-conj-finite-type}}\label{ss:proof-main-result}
Let $\mathcal{A}$ be a cluster algebra  of finite type and $\mathbf{x}$ an arbitrary cluster. We separate the proof by discussing the type of $\mathcal{A}$.

\noindent{\bf Type $\mathbb{A}$}: Assume that $\mathcal{A}$ is of type $\mathbb{A}_n$.  We take $\fS$ to be a disc with $n+3$ boundary marked points and $\fT$ an ideal triangulation of $(\fS,\fM)$ which corresponds to $\mathbf{x}$. Let $\mathbf{m}$ and $\mathbf{n}$ be two cluster monomials such that $d_\mathbf{x}(\mathbf{m})=d_{\mathbf{x}}(\mathbf{n})$. Denote by $\cm$ and $\cn$ the finite multisets corresponding to $\mathbf{m}$ and $\mathbf{n}$, respectively, under the bijection in Theorem \ref{thm:FST-correspondence}. It follows that $\Intv_{\fT}(\cm)=\Intv_{\fT}(\cn)$. By Theorem \ref{t: main theorem}, we conclude that $\cm=\cn$, and hence $\mathbf{m}=\mathbf{n}$ by Theorem \ref{thm:FST-correspondence} again.

\noindent{\bf Type $\mathbb{D}$}: Assume that $\mathcal{A}$ is of type $\mathbb{D}_n$.  We take $\fS$ to be a disc with $n$ boundary marked points and one puncture. Let $\fT$ be a tagged triangulation of $(\fS,\fM,\fP)$  which corresponds to $\mathbf{x}$. Then the result also follows from Theorem \ref{t: main theorem} and Theorem \ref{thm:FST-correspondence} as in type $\mathbb{A}$.

\noindent{\bf Type $\mathbb{B}$ and $\mathbb{C}$}: According to \cite[Theorem 5.10]{FG22}, it suffices to assume that $\mathcal{A}$ is of type $\mathbb{C}_n$. We take $\fS$ to be a disc with $n+1$ boundary marked points and an unmarked boundary component in its interior. According to \cite[Theorem 11]{CL20}, it suffices to consider cluster monomials without involving initial cluster variables.
Then the result follows from Theorem \ref{thm:type-c-intersection-vector} and Proposition \ref{prop:type-c-FST-correspondence}.

\noindent{\bf Exceptional types}: 
For a cluster $\mathbf{x}_s=(x_{1;s},\dots, x_{n;s})$, the $D$-matrix $D_{\mathbf{x}_s}$ of $\mathbf{x}_s$  is the matrix whose $k$-th column vector is the denominator vector of $x_{k;s}$ with respect to an initial cluster $\mathbf{x}$. For a vector $\alpha=(a_1,\dots, a_n)\in \mathbb{N}^n$, we write $\mathbf{x}_s^\alpha=x_{1;s}^{a_1}\cdots x_{n;s}^{a_n}$ for the cluster monomial.
Then we have $d_{\mathbf{x}}(\mathbf{x}_s^\alpha)=D_{\mathbf{x}_s}\alpha^{tr}$, where $\alpha^{tr}$ is the transpose of $\alpha$.

Conjecture \ref{conj:denominator-conj} is false for a cluster algebra $\mathcal{A}$ if and only if there exist two different cluster monomials with the same denominator vector with respect a given cluster $\mathbf{x}$. In particular, one of the following cases occurs:
\begin{itemize}
    \item[Case 1:] There is a cluster $\mathbf{x}_t$ and two vector $\mathbf{u}\neq \mathbf{v}\in \mathbb{N}^n$ such that $d_{\mathbf{x}}(\mathbf{x}_t^\mathbf{u})=d_{\mathbf{x}}(\mathbf{x}_t^\mathbf{v})$. This is equivalent to that $|D_{\mathbf{x}_t}|=0$.
    \item[Case 2:] Suppose that for any cluster $\mathbf{z}$ of $\mathcal{A}$, $|D_{\mathbf{z}}|\neq 0$. Then there exists two different cluster $\mathbf{x}_t$ and $\mathbf{x}_s$ (we may assume that $\mathbf{x}_t\cap \mathbf{x}_s=\{x_{1;t},\dots,x_{r;t}\}$ for some $r<n$), and $\mathbf{u}, \mathbf{v}\in \mathbb{N}^n$ such that $d_{\mathbf{x}}(\mathbf{x}_t^{\mathbf{u}})=d_{\mathbf{x}}(\mathbf{x}_s^{\mathbf{v}})$, where at least one of the last $n-r$ components of $\mathbf{u}$ (resp. $\mathbf{v}$) is non-zero. Equivalently, $D_{\mathbf{x}_t}\mathbf{u}^{tr}=D_{\mathbf{x}_s}\mathbf{v}^{tr}$ and hence  $D_{\mathbf{x}_s}^{-1}D_{\mathbf{x}_t}\mathbf{u}^{tr}=\mathbf{v}^{tr}$. This is further equivalent to there being an $r<l\leq n$ such that the  system of linear inequalities
    \[
    \begin{cases}
        D_{\mathbf{x}_s}^{-1}D_{\mathbf{x}_t}X\geq 0,\\
        X\geq e_l,
    \end{cases}
    \]
     has a solution, where $e_1,\dots, e_n$ is the standard $\mathbb{Z}$-basis of $\mathbb{Z}^n$.
\end{itemize}

Therefore, we may verify Conjecture \ref{conj:denominator-conj} for exceptional types by the following algorithm:

  \begin{enumerate}
      \item[$\triangleright$]Input:  A skew-symmetrizable integer matrix $B\in M_n(\mathbb{Z})$ of finite type.
      \item[Step 1:] Compute the set $\text{Mut}(B)$ of equivalence classes of matrices which can be obtained from $B$ by mutations, which is a finite set, say $\text{Mut}(B)=\{B_1,\dots, B_m\}$.
      \item[Step 2:] For each $B_i\in \text{Mut}(B)$, compute the set of $D$-matrices associated with a cluster pattern of $B_i$, which is a finite set, say $\mathcal{D}_i=\{D_{1;i},\dots, D_{s_i;i}\}$.
      \item[Step 3:] For each $1\leq i\leq m$, compute the determinant $|D_{j;i}|$ for each $D_{j;i}\in \mathcal{D}_i$. If there exist $i$ and $j$ such that $|D_{j;i}|=0$, then Conjecture \ref{conj:denominator-conj} is false. Otherwise,
      \item[Step 4:] For each $1\leq i\leq m$ and for any pair $D_{j;i}, D_{k;i}$, by applying permutations of columns, assume that $D_{j;i}$ and $ D_{k;i}$ have precisely $r$ common columns which are the first $r$ columns.  Solve the following system of  linear inequalities:
      \[
      \begin{bmatrix}
          D_{k;i}^{-1}D_{j;i}\\ E_n
      \end{bmatrix}X\geq \begin{bmatrix}
          0\\ e_l
      \end{bmatrix}
      \]
      for each $r< l\leq n$. 
      \item[$\triangleright$] Conclusion:
      If there are $i,j\neq k,l$ such that the associated system of linear inequalities has a solution, then Conjecture \ref{conj:denominator-conj} is false. Otherwise, Conjecture \ref{conj:denominator-conj} is true for the cluster algebra associated to $B$.
  \end{enumerate}
 The verification and SAGE program for exceptional types can be found in 
  
\url{https://github.com/Changjianfu/denominator-conjecture-exceptional-type}. 
\subsection{Simplification of the algorithm}
We collect certain results which enable us to simplify our verification of  Conjecture \ref{conj:denominator-conj} for cluster algebras of exceptional types. According to \cite{CK08} and \cite{RS20}, it remains to verify Conjecture \ref{conj:denominator-conj} for non-acyclic initial seeds. In particular, Conjecture \ref{conj:denominator-conj} is true for the cluster algebra of type $G_2$, since every seed is acyclic in this case.

 Let $B=(b_{ij})\in M_n(\mathbb{Z})$ be a skew-symmetrizable matrix. If $b_{ik}\geq 0$ for any $1\leq i\leq n$, then the mutation $\mu_k$ in direction $k$ is called a {\em sink mutation}. If $b_{ik}\leq 0$ for any $1\leq i\leq n$, then the mutation $\mu_k$ in direction $k$ is called a {\em source mutation}. Two skew-symmetrizable integer matrices $B_1$ and $B_2$ are {\em sink-source equivalent} if one can obtain $B_2$ from $B_1$ by a finite sequence of sink/source mutations.
The following is a direct consequence of \cite[Theorem 2.4]{RS18}.
\begin{lemma}\label{lem:reduce-to-sink-source-equivalent}
    Let $\mathcal{A}$ be a cluster algebra of finite type. Let $\mathbf{x}$ and $\mathbf{x}'$ be clusters of $\mathcal{A}$ which are related by a sink/source mutation. Then Conjecture \ref{conj:denominator-conj} holds for $\mathcal{A}$ with respect to $\mathbf{x}$ if and only if it holds for $\mathcal{A}$ with respect to $\mathbf{x}'$.
\end{lemma}
According to Lemma \ref{lem:reduce-to-sink-source-equivalent}, it suffices to check Conjecture \ref{conj:denominator-conj} for exchange matrices up to sink-source equivalence.
For example, there are only $5$ non-acyclic exchange matrices of type $F_4$ up to sink-source equivalence, and there are $390$ non-acyclic exchange matrices of type $E_8$ up to sink-source equivalence.

The following lemma reduces the verification of type $E$ to the case $E_8$.
\begin{lemma}\label{lem:sub-cluster-algebra}
    Let $B\in M_n(\mathbb{Z})$ be a skew-symmetrizable matrix and $\mathcal{A}(B)$ the cluster algebra associated with $B$. Denote by $\bar{B}$  the submatrix of $B$ formed by the first $r$ columns of $B$ and $\mathcal{A}(\bar{B})$ the cluster algebras associated with $\bar{B}$. If Conjecture \ref{conj:denominator-conj} is true for $\mathcal{A}(B)$, then Conjecture \ref{conj:denominator-conj} is true for $\mathcal{A}(\bar{B})$.
\end{lemma}
\begin{proof}
 Every cluster monomial of $\mathcal{A}(\bar{B})$ is a cluster monomial of $\mathcal{A}(B)$. Moreover, if $\mathbf{d}$ is the denominator vector of a cluster monomial $z$ of $\mathcal{A}(\bar{B})$, then $\tilde{d}=(\mathbf{d},0,\dots, 0)\in \mathbb{Z}^n$ is the denominator vector of $z$ as a cluster monomial of $\mathcal{A}(B)$ by \cite[Theorem 11]{CL20}. 
\end{proof}

The following lemma greatly reduces the amount of calculation in Step $4$. Namely, it suffices to solve the system of linear inequalities associated with clusters satisfying the condition in Lemma \ref{lem:non-intersection}.
\begin{lemma}\label{lem:non-intersection}
    Let $\mathcal{A}$ be a cluster algebra of rank $n$. Suppose that Conjecture \ref{conj:denominator-conj} is false for $\mathcal{A}$ with respect to an initial cluster $\mathbf{x}=(x_1,\dots, x_n)$. Then there are two different clusters $\mathbf{x}_t=(x_{1;t},\dots, x_{n;t})$ and $\mathbf{x}_s=(x_{1;s},\dots, x_{n;s})$ and two nonzero vectors $\mathbf{u},\mathbf{v}\in \mathbb{N}^n$ satisfy the following conditions: $\mathbf{x}_t\cap \mathbf{x}=\emptyset$, $\mathbf{x}_s\cap \mathbf{x}=\emptyset$, $\mathbf{x}_t\cap \mathbf{x}_s=\emptyset$ and $d_\mathbf{x}(\mathbf{x}_t^\mathbf{u})=d_\mathbf{x}(\mathbf{x}_s^\mathbf{v})$.
\end{lemma}
\begin{proof}
By assumption, there exist two different cluster monomials
$\mathbf{x}_{t_1}^{\mathbf{u}}$ and $\mathbf{x}_{s_1}^{\mathbf{v}}$ with $d_\mathbf{x}(\mathbf{x}_{t_1}^{\mathbf{u}})=d_\mathbf{x}(\mathbf{x}_{s_1}^{\mathbf{v}})$, where $\mathbf{u}=(u_1,\dots,u_n)$ and $\mathbf{v}=(v_1,\dots,v_n)$. We may assume that if $x_k\in \mathbf{x}_{t_1}$ (resp. $\mathbf{x}_{s_1}$), then the multiplicity of $x_k$ in the cluster monomial $\mathbf{x}_{t_1}^\mathbf{u}$ (resp. $\mathbf{x}_{s_1}^\mathbf{v}$) is zero. Indeed, assume that $x_k\in \mathbf{x}_{t_1}$, say $x_k=x_{1;t_1}$, and $u_1>0$. Denote by $\mathbf{u}'=\mathbf{u}-u_1e_1$. We have $d_\mathbf{x}(\mathbf{x}_{t_1}^{\mathbf{u}})=-u_1e_1+d_\mathbf{x}(\mathbf{x}_{t_1}^{\mathbf{u}'})$. By Theorem \cite[Theorem 11]{CL20}
and $d_\mathbf{x}(\mathbf{x}_{t_1}^\mathbf{u})=d_\mathbf{x}(\mathbf{x}_{s_1}^\mathbf{v})$, we conclude that $x_k\in \mathbf{x}_{s_1}$. Without loss of generality, we may assume that $x_k=x_{1;s_1}$. If follows that $u_1=v_1$. Set $\mathbf{v}'=\mathbf{v}-v_1e_1$. Then, we can replace $\mathbf{x}_{t_1}^{\mathbf{u}}$ by $\mathbf{x}_{t_1}^{\mathbf{u}'}$ and $\mathbf{x}_{s_1}^{\mathbf{v}}$ by $\mathbf{x}_{s_1}^{\mathbf{v}'}$.

 Assume that $x_k\in \mathbf{x}_{t_1}$. By assumption, the multiplicity of $x_k$ in $\mathbf{x}_{t_1}^\mathbf{u}$ is zero. We may mutate $\mathbf{x}_{t_1}$ at $x_k$ to obtain a new cluster $\mathbf{x}_{t_2}=\mathbf{x}_{t_1}\backslash\{x_k\}\cup \{x_k'\}$. According to \cite[Theorem 4.22]{FuGy24}, we have $(x_k||x_k')=1$, where $(x_k||x_k')$ is the $f$-compatibility degree between $x_k$ and $x_k'$. It follows that $x_k'\not\in \mathbf{x}$ by \cite[Theorem 4.18]{FuGy24}. Consequently, $\mathbf{x}_{t_2}\cap \mathbf{x}=(\mathbf{x}_{t_1}\cap \mathbf{x})\backslash \{x_k\}$ and  $\mathbf{x}_{t_1}^{\mathbf{u}}=\mathbf{x}_{t_2}^{\mathbf{u}}$. By repeating the process, we obtain a cluster $\mathbf{x}_{t_3}$ such that $\mathbf{x}_{t_3}^{\mathbf{u}}=\mathbf{x}_{t_1}^{\mathbf{u}}$ and $\mathbf{x}_{t_3}\cap \mathbf{x}=\emptyset$. Similarly, we have a cluster $\mathbf{x}_{s_3}$ such that $\mathbf{x}_{s_3}^{\mathbf{v}}=\mathbf{x}_{s_1}^{\mathbf{v}}$ and $\mathbf{x}_{s_3}\cap \mathbf{x}=\emptyset$.

 If $\mathbf{x}_{t_3}\cap \mathbf{x}_{s_3}\neq \emptyset$. Without loss of generality, we may assume that $x_{1;t_3}=x_{1;s_3}$ and $u_1\geq 0$, $v_1=0$. Indeed, assume that $u_1\geq v_1$ and set $\mathbf{u}'=\mathbf{u}-v_1e_1$, $\mathbf{v}'=\mathbf{v}-v_1e_1$.  We may then replace $\mathbf{x}_{t_3}^{\mathbf{u}}$ by $\mathbf{x}_{t_3}^{\mathbf{u}'}$ and $\mathbf{x}_{s_3}^{\mathbf{v}}$ by $\mathbf{x}_{s_3}^{\mathbf{v}'}$, respectively. By considering the mutation of $\mathbf{x}_{s_3}$ at $x_{1;s_3}$, we obtain a new cluster $\mathbf{x}_{s_4}=\mathbf{x}_{s_3}\backslash\{x_{1;s_3}\}\cup\{x_{1;s_3}'\}$. Again, by \cite[Theorem 4.18 \& 4.22]{FuGy24}, we conclude that $x_{1;s_3}'\not\in \mathbf{x}_{t_3}$.
It follows that  $|\mathbf{x}_{t_3}\cap\mathbf{x}_{s_4}|=|\mathbf{x}_{t_3}\cap\mathbf{x}_{s_3}|-1$. Note that $\mathbf{x}_{s_4}^{\mathbf{v}}=\mathbf{x}_{s_3}^\mathbf{v}$. Repeating the process, we obtain the desired cluster monomials. This completes the proof.
\end{proof}

\subsection{A consequence in cluster-tilted algebras}\label{ss: cluster-tilted-algebra}
Let $Q$ be a finite acyclic quiver and $\cc_Q$ the associated cluster category \cite{BMRRT} (cf. \cite{CCS2006} for $A_n$ cases). We denote by $\Sigma$ the suspension functor of $\cc_Q$.
Recall that an object $M\in \cc_Q$ is {\em rigid} provided that $\Ext^1_{\cc_Q}(M,M)=0$.
An object $T\in \cc_Q$ is a {\em cluster-tilting} object if $T$ is rigid and satisfies tha $\Ext^1_{\cc_Q}(T,X)=0$  implies that $X\in \operatorname{add} T$, where $\operatorname{add} T$ is the full subcategory consisting of objects  which are finite direct sums of direct summands of $T$. The endomorphism algebra $\End_{\cc_Q}(T)$ of a basic cluster-tilting object is called a {\em cluster-tilted algebra} \cite{BMR2007}. A cluster-tilted algebra is of finite representation type if and only if the underlying graph of $Q$ is of type $\mathbb{A},\D$ or $\E$.

Denote by $Q_0=\{1,\dots, n\}$ the vertex set of $Q$. Let  $\ma(Q)$ be the cluster algebra associated to $Q$.
According to~\cite{CK06}, there is a bijection $\varphi$ between the set of cluster variables of $\ma(Q)$ and the set of isomorphism classes of indecomposable rigid objects of $\mc_Q$, which induces a bijection between the set of clusters of $\ma(Q)$ and the set of isomorphism classes of basic cluster-tilting objects of $\mc_Q$. Fix a cluster $\mathbf{y}=\{y_1,\dots, y_n\}$ of $\ma(Q)$, denote by $\Sigma T=\bigoplus_{i=1}^n(\Sigma T_i)\in \mc_Q$ the corresponding basic cluster-tilting object, where $\Sigma T_1,\dots, \Sigma T_n$ are indecomposable rigid objects corresponding to $y_1,\dots, y_n$ respectively. Let $\Gamma=\End_{\mc_Q}(T)$ be the cluster-tilted algebra and $\mod \Gamma$ the category of finitely generated right $\Gamma$-modules. A module $M\in \mod \Gamma$ is {\em $\tau$-rigid} if $\Hom_{\Gamma}(M,\tau M)=0$, where $\tau$ is the Auslander-Reiten translation. A pair $(M,P)$ is a {\em $\tau$-rigid pair} if $M$ is $\tau$-rigid and $P$ is a finitely generated projective $\Gamma$-module such that $\Hom_\Gamma(P,M)=0$.
By \cite[Theorem 4.1]{AIR14}, the bijection $\varphi$ induces a bijection between the set of isomorphism classes of $\tau$-rigid pairs of $\Gamma$-modules and the set of cluster monomials of $\ma(Q)$. By abuse of notation, we still denote the induced bijection by $\varphi$. Each $\tau$-rigid module $M$ can be viewed as a $\tau$-rigid pair $(M,0)$. Hence for each $\tau$-rigid module $M$, $\varphi(M)$ is a cluster monomial of $\ma(Q)$.
If $\Gamma$ is of finite representation type,  we have $d_\mathbf{y}(\varphi(M))=\dimv M$  for each $\tau$-rigid $\Gamma$-module  $M$ by \cite{CCS06}.

 As a direct consequence of  Theorem~\ref{thm:denominator-conj-finite-type}, we have
\begin{thm}\label{thm:dimension-vector} 
Let $\Gamma$ be a cluster-tilted algebra of finite representation type. Then different $\tau$-rigid $\Gamma$-modules have different dimension vectors.
    \end{thm}

\begin{remark}
Let $\Gamma$ be a cluster-tilted algebra of finite representation type.
It has been proved in \cite{GP12,R11} that different indecomposable $\Gamma$-modules have different dimension vectors. Note that every indecomposable $\Gamma$-module is $\tau$-rigid in this case. Hence, Theorem \ref{thm:dimension-vector} is a non-trivial generalization of the corresponding result of \cite{GP12, R11}.
\end{remark}


\bibliographystyle{acm}
\bibliography{myref}

\end{document}